%% file: LSRGB.tex
\documentclass[twocolumn,10pt]{asme2e}
\pdfoutput=1

\usepackage{graphicx}
\usepackage{float}
\usepackage{amsmath}
\usepackage{amsfonts}
\usepackage{amssymb}
\usepackage[thinspace, thinqspace,amssymb]{SIunits}
\usepackage{subfigure}
\usepackage{placeins}
\usepackage{acronym}
\usepackage{pgfplots}
\usepackage{tikz}
\usetikzlibrary{positioning}

\title{Probabilistic LCF Risk Evaluation of a Turbine Vane by Combined Size Effect and Notch Support Modeling}

\author{Lucas M\"ade
    \affiliation{Gas Turbine Department of Materials and Technology\\
		Siemens AG\\
   Berlin, 10553\\
    Germany\\
    Email: lucas.maede@siemens.com\vspace{5mm}
    }\\
    {\tensfb Sebastian Schmitz}
    \affiliation{Gas Turbine Department of Materials and Technology\\
		Siemens AG\\
       Berlin, 10553\\
    Germany\\
    Email: schmitz.sebastian@siemens.com\vspace{5mm}
    }\\
    {\tensfb Georg Rollmann}
    \affiliation{Gas Turbine Department of Materials and Technology\\
		Siemens AG\\
    M\"ulheim an der Ruhr, Nordrhein-Westfalen, 45473\\
    Germany\\
    Email: georg.rollmann@siemens.com
    }
}

\author{Hanno Gottschalk
\affiliation{Faculty of Mathematics and Natural Science\\
    Bergische Universit\"at Wuppertal\\
    Wuppertal, Nordrhein-Westfalen, 42119\\
    Germany\\
    Email: hanno.gottschalk@uni-wuppertal.de\vspace{5mm}
    }\\
       {\tensfb Tilmann Beck}
    \affiliation{Institute of Materials Science and Engineering\\
    Technische Universit\"at Kaiserslautern\\
      Kaiserslautern, Rheinland-Pfalz, 67653\\
    Germany\\
    Email: beck@mv.uni-kl.de
    }
}

\begin{document}

\maketitle    

\begin{abstract}
{\it A probabilistic risk assessment for low cycle fatigue (LCF) based on the so-called size effect has been applied on gas-turbine design in recent years. In contrast, notch support modeling for LCF which intends to consider the change in stress below the surface of critical LCF regions is known and applied for decades. Turbomachinery components often show  sharp stress gradients and very localized critical regions for LCF crack initiations so that a life prediction should also consider notch and size effects. The basic concept of a combined probabilistic model that includes both, size effect and notch support, is presented. In many cases it can improve LCF life predictions significantly, in particular compared to \textit{E-N} curve predictions of standard specimens where no notch support and size effect is considered. Here, an application of such a combined model is shown for a turbine vane. 
}

\end{abstract}

\begin{nomenclature}
\entry{LCF}{Low cycle fatigue}
\entry{\textit{E-N}~curve}{Curve of strain amplitude vs. crack initiation life}
\entry{$CMB$}{Coffin-Manson-Basquin (model)}
\entry{$\sigma_f^\prime$}{Fatigue strength coefficient}
\entry{$b$}{Fatigue strength exponent}
\entry{$\epsilon_f^\prime$}{Fatigue ductility coefficient}
\entry{$c$}{Fatigue ductility exponent}
\entry{$E$}{Cyclic Young's modulus}
\entry{$N_i$}{Load cycles until crack initiation}
\entry{$m$}{Weibull shape parameter}
\entry{$\eta$}{Weibull scale parameter}
\entry{$\chi$}{Normalized gradient of equivalent elastic stress}
\entry{$\chi_T$}{Normalized gradient of temperature}
\entry{$\sigma_e$}{Equivalent elastic stress}
\entry{$NSP$}{Notch support parameters}
\entry{$\Omega$}{3D component domain}
\entry{$\partial\Omega$}{2D domain surface}
\entry{MLE}{Maximum likelihood estimation}
\entry{FEA}{Finite Element Analysis}
\entry{TMF}{Thermo mechanical fatigue}
\entry{TBC}{Thermal barrier coating}
\entry{BC}{Bond coat}
\end{nomenclature}

\input{Introduction}

\input{Section1}

\input{Section2}
\input{Conclusion}
\begin{acknowledgment}
We wish to thank the gas turbine technology department of the Siemens AG for stimulating discussions and many helpful suggestions.
\end{acknowledgment}

\bibliographystyle{asmems4}



\end{document}

%% file: Introduction.tex
\section*{INTRODUCTION}
It is well known that the number of cycles until initiation of a fatigue crack in Ni-based superalloys is subjected to considerable statistical scatter, see e.g. \cite{Vormwald} for a discussion of statical scatter in fatigue experiments.
The design of gas turbines and their safe and reliable operation therefore requires a methodology that is capable to accurately quantify risk levels for \textit{low cycle fatigue} (LCF) crack initiation, crack growth and ultimate failure. Traditional deterministic design approaches however predict absolute safety below a specified 'safe' number of service cycles and failure just above it. Such life prediction models give clear answers on one hand but on the other hand they do not provide an adequate description of the real world. At the same time, such a binary description narrows the business options of gas turbine power plant operators and service providers, where in some instances taking a non safety relevant economic risk in exchange for an even bigger economic opportunity might be a rational way of decision making. Of course, such decisions need to be underpinned by a proper risk assessment.
Responding to this need, in \cite{Schmitz_Seibel,ASME2013Paper}, a probabilistic model, based on a local Weibull hazard density approach is being used for LCF crack initiation prediction. It inherently considers the statistical size effect and the inhomogeneity of surface stress and has been applied to different gas turbine components, such as blades and compressor discs \cite{LCF7Blade}. Further probabilistic models for LCF have been proposed by other authors, e.g. \cite{Fedelich,HertelVormwald,Beretta}. Okeyoyin \textit{et al.} used a probabilistic framework for computation of fatigue notch factors for high cycle fatigue (HCF) based on a random distribution of failure inducing defects in the volume of material \cite{Okeyoyin}.
While Hertel \textit{et al.} \cite{HertelVormwald} take into account notch support factors, this is not the case for the other papers. The modes \cite{HertelVormwald,Beretta} however are based on the Paris-Erdogan law of crack growth and an initial flaw size distribution in the sense of strength-probability-time (SPT)-diagrams in ceramics. Also, the model \cite{HertelVormwald} has been set up and validated for several steels and not for superalloys. So these models considerably differ from the local probabilistic model for LCF proposed and validated in \cite{Schmitz_Seibel,ASME2013Paper}. However, notch support factors were not included beyond the statistical size effect.
This gap is here closed for the first time where the notch support effect arising from stress gradients in the volume combined with the size effect is included in the probabilistic framework and experimental validation work is provided.
Here the crack initiation prediction prediction of a turbine vane is presented as a use case of the probabilistic model with notch support implementation. Section~\ref{sec:LCF} provides a brief repetition of the key steps to understanding the local probabilistic model with a focus at the notch support mechanism in Subsection~\ref{sec:LCF_NSE}. The following Subsection \ref{subsec:Model_Validation} outlines methodology and results of the notch support model validation for \unit{850}{\degree\Celsius}. A crack initiation prediction of a turbine vane is discussed in Section~\ref{sec:LifePredictVane}, emphasizing the differences of the results when neglecting (Subsections~\ref{sec:TBV_noNSP}) and enabling (Subsection~\ref{sec:TBV_NSP}) notch support in the model.



%% file: Section1.tex
\section{LOW CYCLE FATIGUE}\label{sec:LCF}
Failure due to strain driven LCF surface crack initiation plays an important role for highly loaded engineering parts, such as turbine components made of Ni-based super alloys. Since these alloys feature high yield strength and relatively low ductility, the initiation of surface cracks of a critical length could lead to rapid failure of the component under cyclic load.

The empirical relationship between maximum load (strain amplitude) at the component and the number of load cycles to failure $N_i$, the \textit{E-N}~curve (or Woehler-curve), is the basis for the local probabilistic model discussed in this article. Several deterministic models for Woehler-curves exist. A well known is the \textit{Coffin-Manson-Basquin} (CMB) equation describing Woehler-curves in strain controlled fatigue \cite{Vormwald,Basquin,Coffin,Harders_Roesler},
\begin{align}
&\epsilon_a=\frac{\sigma_f^\prime}{E}\left(2N_i\right)^{b}+\epsilon_f^\prime\left(2N_i\right)^{c}.
\label{eq:LCF_CMB_Gl}
\end{align}
The parameters $\sigma_f,\ b$ (fatigue strength coefficient/exponent) and $\epsilon_f,\ c$ (fatigue ductility coefficient/exponent) are material parameters obtained by fitting test data while $E$ is the cyclic Young's modulus. The first summand at the r.h.s. of Eqn.~\eqref{eq:LCF_CMB_Gl} describes \textit{E-N}~curves dominated by elastic strain, while the second summand at the r.h.s accounts for dominating plastic strain. In the conventional safe life approach of crack initiation prediction the CMB model is used in combination with safety factors to account for natural scatter, size effects and additional effects influencing the LCF mechanism. The probabilistic approach to crack initiation prediction introduced in \cite{Schmitz_Seibel} is explicitly accounting for the first two of these influences. For more reliable and realistic predictions, the \textit{notch support effect} is additionally implemented in this model.

\subsection{Local Probabilistic Model for LCF}\label{sec:LCF_locPM}
The motivation of using a probabilistic model for crack initiation prediction is the need to account for the scatter which is generally observed in fatigue events and the statistical size effect. The latter accounts for the different probabilities of crack initiation in bodies of equal shape but different size when subjected to equal stress.

In a probabilistic framework for crack initiation prediction, the number of cycles to crack initiation $N_i$ is regarded as a random variable whose statistics is characterized by the cumulative distribution function $F_{N_i}(n)$ and probability density function $f_{N_i}(n)$. Here the local approach from \cite{ASME2013Paper} is taken up again. It assumes all members of $\{(N_i)_j\}_{j=1...k}$ for the subsets of an arbitrary partition $\{A_j\}_{j=1...k}$ of the component surface $\partial\Omega$, to be independent from each other because initial cracks only cover the range of few grains. The concept of the \textit{hazard rate} $h(n)$ was chosen to quantify the risk for crack initiation in every surface patch $A_j$ individually as its property of additivity for stochastically independent variables allows partitioning of the risk analysis of the body's surface. That is a requirement of the local approach for crack initiation prediction at the entire surface of an arbitrarily shaped body. The hazard rate is calculated by
\begin{align}
h(n)=\lim\limits_{\Delta n \to 0} \frac{P(n<N_i<n+\Delta n\mid N_i>n)}{\Delta n}=\frac{f_{N_i}(n)}{1-F_{N_i}(n)}.
\end{align}
Thus $h(n)\cdot \Delta n$ is the probability of crack initiation within the cycle $n+\Delta n$ where $\Delta n$ is the cycle increment \cite{Escobar_Meeker}. The hazard rate for the first crack initiation at the entire surface $\partial\Omega$ is the sum of those for all $A_j$, $h=\sum_{j=1}^{k}{h_j}$ since one assumes the number of load cycles until crack initiation $(N_i)_j$ to be stochastically independent in every $A_j$. The second assumption states that crack initiation risk is a functional of only local strain and temperature fields $\epsilon_a(\mathbf{x})$ and $T(\mathbf{x})$ since no long range order phenomena occur in the continuum mechanics model for polycrystalline materials. That is why $h(n)$ is the surface integral of $\rho(n;\epsilon_a,T)$ the in the limit of an infinitesimal fine partition $\{A_i\}_i$ with $\rho(n;\epsilon_a,T)$ being a local functional of strain and temperature field at the surface $\partial\Omega$.
\begin{align}
h(n)=\int_{\partial\Omega}{\rho(n;\epsilon_a(\mathbf{x}),T(\mathbf{x}))\,dA}\label{eq:SSE_HazDenInt}
\end{align}
The \textit{cumulative hazard function} $H(n)$ is defined as 
\begin{equation}
H(n)=\int_{0}^{n}{h(t)\, dt}.\label{eq:SSE_cumHazFunc}
\end{equation}
$H(n)$ and $h(n)$ fulfill the relations
\begin{align}
h(n)&=\frac{f_{N_i}(n)}{1-F_{N_i}(n)},\quad F_{N_i}(n)=1-e^{-H(n)}.\label{eq:SSE_cumHazDistFunc}
\end{align}
In \cite{Schmitz_Seibel} the number of cycles to crack initiation $N_i$ is assumed to be Weibull distributed with the cumulative distribution function
\begin{align}
&F_{N_i}(n)=1-\mathrm{exp}\left[-\left(\frac{n}{\eta}\right)^m\right],\label{eq:locPM_newCDF}
\end{align}
where $m$ is the Weibull shape parameter and $\eta$ the Weibull scale parameter. While $\eta$ determines the position of the distribution in the domain, $m$ influences the shape of the distribution (broad or peaked) and thus the expected scatter of events. The whole concept of the local approach leads to an integration formula for $\eta$ which adds up the risk for crack initiation along the examined surface.
Using equations~\eqref{eq:SSE_cumHazDistFunc} for the Weibull distribution one finds
\begin{align}
\int_{\partial\Omega}{\rho(n;\epsilon_a(\mathbf{x}),T(\mathbf{x}))\,dA}=\frac{m}{\eta}\cdot \left(\frac{n}{\eta}\right)^{m-1}.
\end{align}
The local approach effectively states that all distributions for $(N_i)_j$ scale individually dependent on the load state. Consequently the integrand in Eqn.~\eqref{eq:SSE_HazDenInt} is
\begin{align}
\rho(n;\epsilon_a,T)=\frac{m}{N_{i_\mathrm{det}}(\epsilon_a,T)}\left(\frac{n}{N_{i_\mathrm{det}}(\epsilon_a,T)}\right)^{m-1}.\label{eq:locPM_rho}
\end{align}
There $N_{i_\mathrm{det}}(\epsilon_a(\mathbf{x}),T(\mathbf{x}))$ is the deterministic number of life cycles at every point $\mathbf{x}$ of the body's surface $\partial\Omega$.

The integrand Eqn.~\eqref{eq:locPM_rho} allows independent integrations over surface and time, which are necessary to receive the cumulative hazard function $H(n)$ according to equations~\eqref{eq:SSE_HazDenInt} and \eqref{eq:SSE_cumHazFunc}. The cumulative hazard function is then found to have the formula
\begin{align}
H(n)=n^m\cdot\int_{\partial\Omega}{\frac{1}{N_{i_{\mathrm{det}}}^m}}\ dA,\label{eq:cumHazFunc}
\end{align}
where the remaining integrand $\left(1/N_{i_{\mathrm{det}}}(\mathbf{x})\right)^m$ is defined as \textit{hazard density}. 

In the following Section \ref{sec:LifePredictVane} plots of the hazard density field at the geometry are shown in order to visualize the risk at the components surface. This is preferable compared to plots of $N_{i_{\mathrm{det}}}$ as one can add up values of $\rho(n;\epsilon_a(\mathbf{x}),T(\mathbf{x}))$ from arbitrary surface spots to receive the overall hazard density for the combined surface. This is a convenient way to directly assess the criticality of different surface subsets. From equations~\eqref{eq:SSE_cumHazDistFunc}, \eqref{eq:locPM_newCDF} and \eqref{eq:cumHazFunc} one further derives a formula for the Weibull scale parameter: 
\begin{align}
\eta=\left(\int_{\partial\Omega}{\frac{1}{N_{i_\mathrm{det}}^m}dA}\right)^{-1/m}.\label{eq:locPM_ScaleSurfInt}
\end{align}
The deterministic life $N_{i_\mathrm{det}}(\mathbf{x})$ of one surface patch is determined by numerically solving
\begin{align}
\epsilon_a(\mathbf{x})=\frac{\sigma_f^\prime}{E}\left(2N_{i_\mathrm{det}}(\mathbf{x})\right)^{b}+\epsilon_f^\prime\left(2N_{i_\mathrm{det}}(\mathbf{x})\right)^{c}.\label{eq:locPM_locCMB}
\end{align}
The computational realization of this model is a tool that uses Finite Element Models (FEA) as input and computes $N_{i_\mathrm{det}}$ at all integration points. The parameters $\sigma_f^\prime,\ \epsilon_f^\prime,\ b,\ c$ for Eqn.~\eqref{eq:locPM_locCMB} are now valid for one small surface patch and are thus independent of the investigated component geometry, given that the mesh of the FEA input is sufficiently small. Hence one can also interpret them as material parameters. They are simultaneously derived with $m$ from maximum likelihood fits of specimen test data. Shape $m$ and scale $\eta$ entirely define the distribution function in Eqn.~\eqref{eq:locPM_newCDF} from which the \unit{50}{\%}-quantile is regarded as the probabilistic average life until crack initiation. By integrating $1/N_{i_\mathrm{det}}^m(\mathbf{x})$ over the entire surface in Eqn.~\eqref{eq:locPM_ScaleSurfInt}, the presented local probabilistic model inherently incorporates the statistical size effect and accounts for material scatter through the shape parameter $m$.


\subsection{Notch Support Effect}\label{sec:LCF_NSE}
Components with notches or other inhomogeneous geometry features exhibit domains of concentrated stress at the respective location when subjected to a load. Geometry induced stress concentration leads to inhomogeneous stress fields in the affected domain while the highest values are usually occurring at the surface. Whereas domains near the surface quickly reach yield strength and are therefore plastically strained, domains further inside the body still support the structure since they experience smaller stresses and therefore impede failure. That is why the crack initiation life of parts exhibiting spatially inhomogeneous stress fields under cyclic load is higher than predicted by the  CMB equation for the maximum occurring strain $\epsilon_a$. Siebel \textit{et al.} have approached a quantification of this phenomenon, known as \textit{notch support effect}, with a support number $n_{NS}$
\begin{equation}
n_{NS}=\frac{\text{observed fatigue strength}}{\text{expected fatigue strength}}=\frac{\sigma_\mathrm{notched/obs}}{\sigma_\mathrm{homogeneous}}.\label{eq:NSE_notchsupportnu}
\end{equation}
They considered $n_{NS}$ to be directly proportional to the stress gradient in the loaded component \cite{Siebel}. This is a well justified approach since quickly abating loads (high concentration) imply larger low stress domains to support the structure. Hence, a stress gradient based support factor $n_\chi$ is also used to consider the notch support effect in the use case described in this paper, where $n_\chi$ is dependent on $\chi(\mathbf{x})$ and material dependent \textit{notch support parameters}\footnote{See Fig.~10.36 on page 378 in \cite{Harders_Roesler} for relationship between $n_\chi$ and $\chi$} $A$ and $k$. They are simultaneously derived with the CMB parameters from LCF test data as described in Section \ref{subsec:Model_Validation}.
\begin{align}
\chi(\mathbf{x})=\frac{1}{\sigma_e(\mathbf{x})}\,\nabla\sigma_e(\mathbf{x})\ \ \mathrm{with}\ \ \mathbf{x}\in\partial\Omega\label{eq:NSE_Chicalc}
\end{align}
is the derivative of the elastic von Mises stress $\sigma_e$ normalized with its surface value. Note that $\sigma_e$ is a scalar field in the whole domain $\Omega$ so that the gradient is well defined at the surface $\partial\Omega$.\\
Since $n_\chi$ can be seen as the strain equivalent to Eqn.~\eqref{eq:NSE_notchsupportnu}, it is combined with the CMB equation to Eqn.~\eqref{eq:locPM_locCMB_andNS} instead of Eqn.~\eqref{eq:locPM_locCMB} for $N_{i_\mathrm{det}}(\mathbf{x})$ computation. This shifts the \textit{E-N}~curve to higher life because $n_\chi\geq 1$. Then, Weibull scale $\eta$, computed by integrating $1/N_{i_\mathrm{det}}^m(\mathbf{x})$ over the surface, and shape $m$ define a distribution $F_{N_i}(n)$ for load cycles $n$ until crack initiation that accounts simultaneously for size effect and notch support effect.\\
Thus, the probabilistic model for LCF with combined size effect and notch support is given by the following Weibull approach:
\begin{align}
F_{N_i}(n)&=1-\mathrm{exp}\left[-\left(\frac{n}{\eta}\right)^m\right]\label{eq:DistrFunct},\\
\eta&=\left(\int_{\partial\Omega}{\frac{1}{N_{i_\mathrm{det}}^m}dA}\right)^{-1/m}\label{eq:etaInt_HazDen},\\
\frac{\epsilon_a(\mathbf{x})}{n_\chi(\mathbf{x})}&=\frac{\sigma_f^\prime}{E}\left(2N_{i_\mathrm{det}}(\mathbf{x})\right)^{b}+\epsilon_f^\prime\left(2N_{i_\mathrm{det}}(\mathbf{x})\right)^{c}.\label{eq:locPM_locCMB_andNS}
\end{align}
Note that besides of the notch support effect, the statistical size effect is also playing an important role in the LCF life of irregularly shaped components because critical stresses usually occur in confined domains which are small compared to the entire component.


\subsection{Calibration of Notch Support Parameters and Model Validation}\label{subsec:Model_Validation}
As mentioned in the previous subsection, the CMB-, notch support- and shape parameters for the combined local probabilistic model are estimated from material test data. The principle of the model calibration and validation procedure is shortly described here.

\begin{figure}[htbp]
		\centering
		\includegraphics[trim=0cm 0cm 1cm 2cm,clip,width=0.45\textwidth]{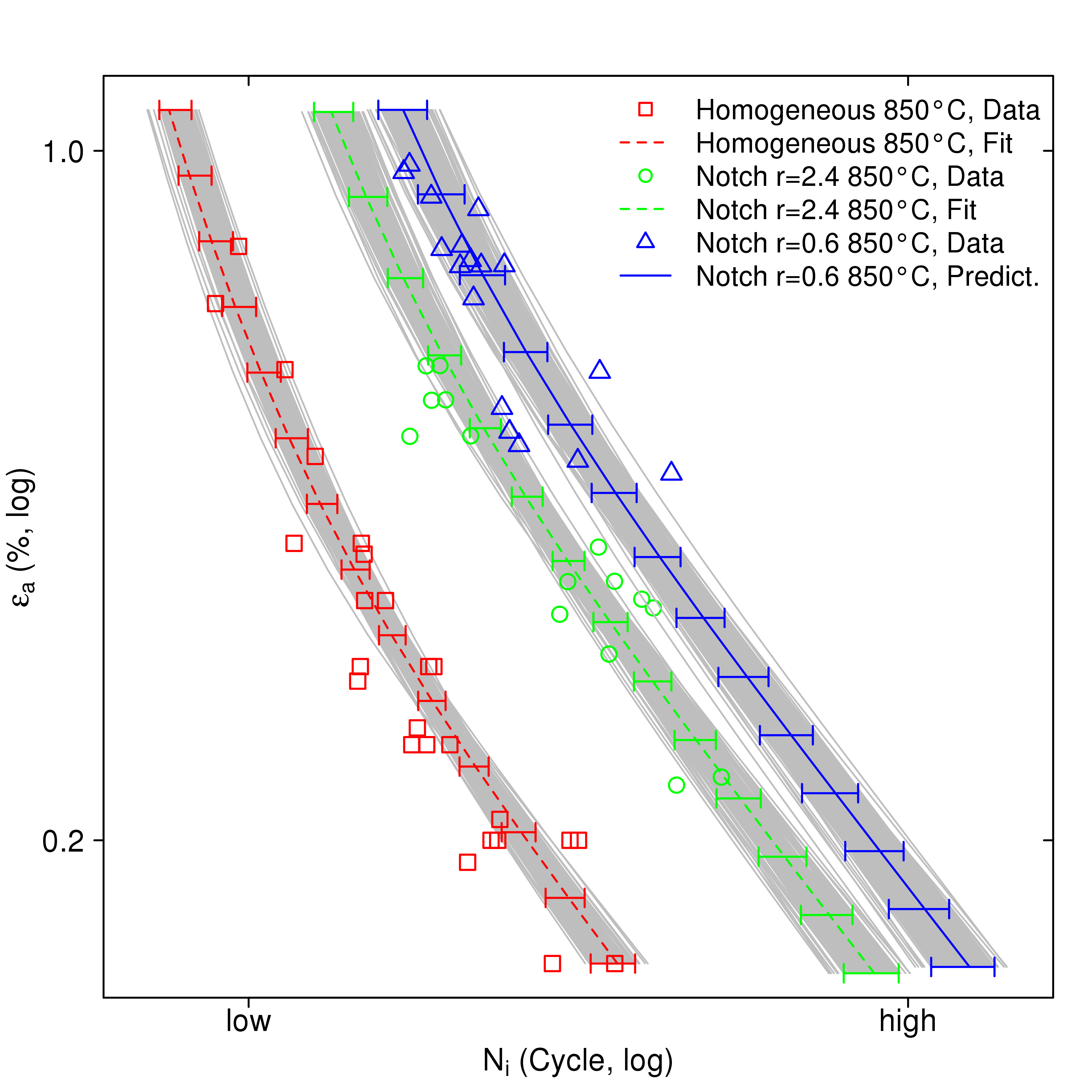}
		\caption{\textit{E-N} CURVES AND TEST DATA OF HOMOGENOUS AND NOTHCED SPECIMEN: FIT AND PREDICTION}
		\label{fig:Fit_PredCurve}
		\vspace{-1\baselineskip}
\end{figure}

In order to calibrate the notch support model, LCF-test data of a specimen with homogeneous, cylindrical gauge section (red squares in Fig.~\ref{fig:Fit_PredCurve} \cite{Schmitz_Seibel}) is simultaneously fitted with test data of a specimen with a circumferential notch of radius \unit{2.4}{mm} (green circles in Fig.~\ref{fig:Fit_PredCurve} \cite{Seibel_Diss}). Both sample types are made of polycrystalline cast RENE 80 and tested in strain controlled LCF at \unit{850}{\celsius}. The local probabilistic model extended with notch support in Subsection \ref{sec:LCF_NSE} requires the parameter set
\begin{equation}
\mathbf{\theta}=\left(\sigma_f,b,\epsilon_f,c,A,k,m\right)^T
\end{equation}
which is determined by maximum likelihood estimation (MLE). Apart from parameters $A$ and $k$ in $\mathbf{\theta}$, the fitting procedure follows \cite{Schmitz_Seibel} from this point on. Determining the minimum of the negative log-likelihood function is achieved with Nelder-Mead optimization. The parameter estimate $\mathbf{\hat{\theta}}$ is then used to compute the median \textit{E-N}~curve of another notched specimen with notch radius \unit{0.6}{mm} (solid blue line). All test data points and \textit{E-N}~curves are shown in Fig.~\ref{fig:Fit_PredCurve}.

The error bars in Fig.~\ref{fig:Fit_PredCurve} indicate the \unit{92.5}{\%}-confidence interval of the respective median life. Uncertainties in parameter estimation are computed by parametric bootstrapping. 2000 bootstrap sample\footnote{The minimum number of bootstrap samples according to \cite{Escobar_Meeker} is chosen for feasible computation times.} sets are generated from the original crack initiation life distributions and fitted with the mentioned MLE procedure. 200 of the resulting median \textit{E-N}~curves are plotted in grey for each specimen geometry. The confidence intervals here represent the \unit{92.5}{\%} percentile range of the uncertainty in curve prediction. Hence they do not cover the observed residual scatter which is e.g. in this case of superalloys also dependent on the grain orientation and the location of the initial cracks compared to the probing tips of the extensometer in the LCF experiment. The calibration curves (dashed) show that the observed test data is well described by the current notch support model. Additionally, the solid blue \textit{E-N}~curve, a pure prediction, shows a good validation for another notch specimen data set (blue triangles) not used for calibration. These findings verify the appropriability of the $\chi$-approach for the available test data. This motivates the application of the notch support model to a turbine vane discussed in Subsection \ref{sec:TBV_NSP}. 



%% file: Section2.tex
\section{PROBABILISTIC LIFE PREDICTION FOR A TURBINE VANE UNDER THERMOMECHANICAL LOADS}\label{sec:LifePredictVane}
\begin{figure}[htbp]
		\centering
		\subfigure[von Mises stress]{
		\begin{tikzpicture}
			\node[anchor=south west,inner sep=0] (StressMap) at (0,0) {\includegraphics[trim=13cm 10mm 13.6cm 2.5cm,clip,width=0.235\textwidth]{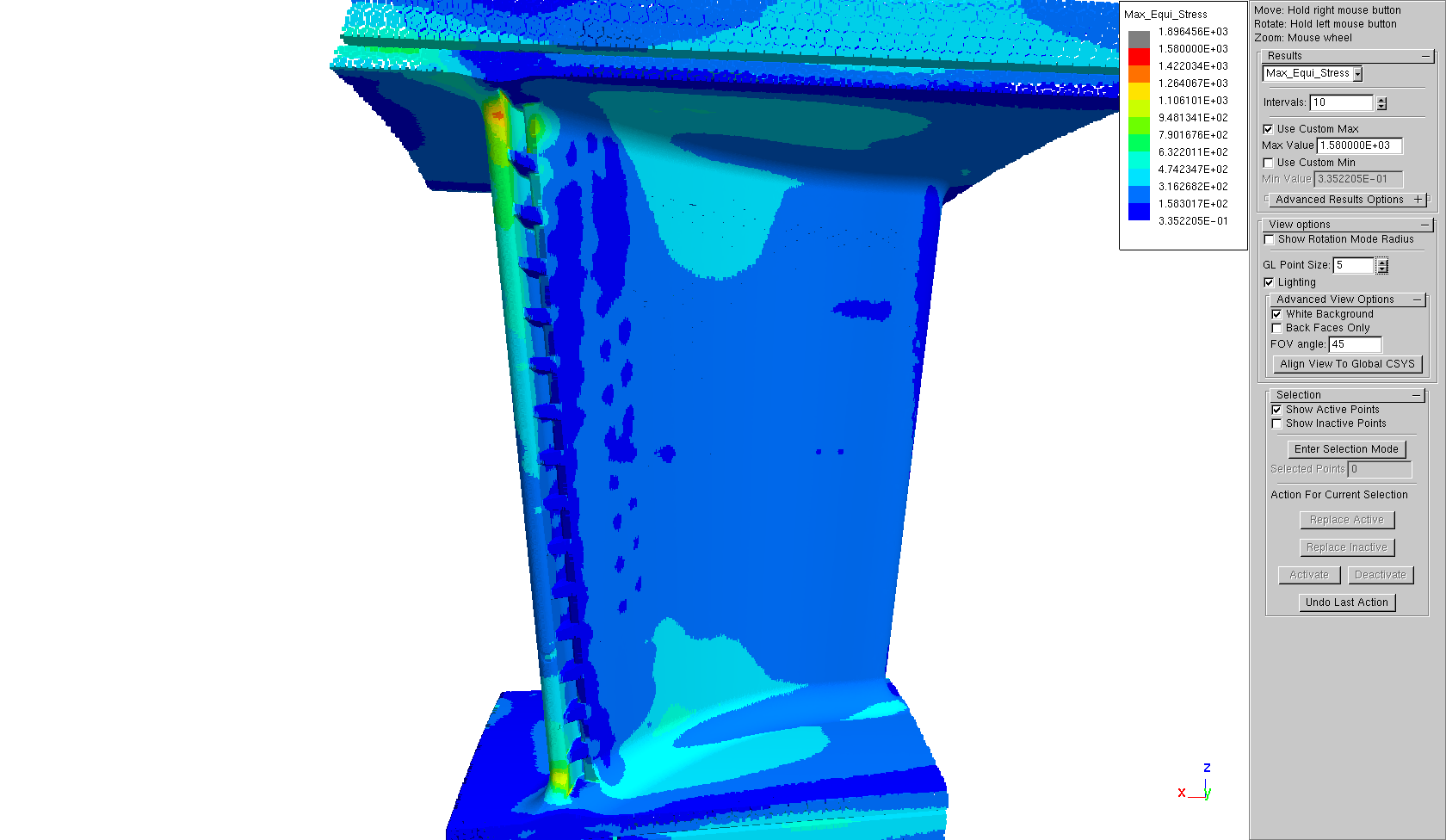}};
		\begin{scope}[x={(StressMap.south east)},y={(StressMap.north west)}]
		\node[anchor=south west,inner sep=0] at (0.95,0) {\includegraphics[trim=77.5cm 35.6cm 12cm 1.9cm,clip,width=0.007\textwidth]{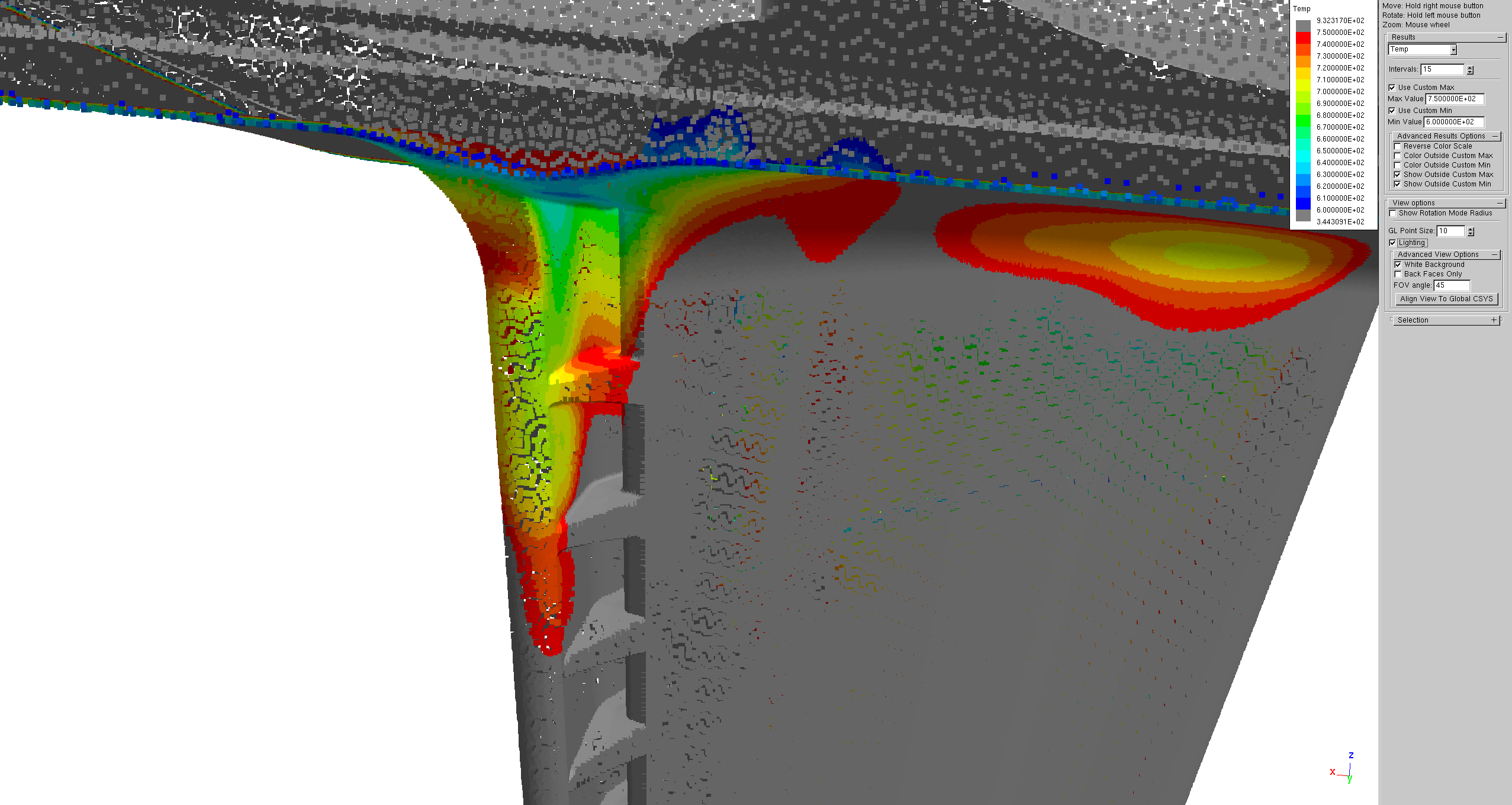}};
		\node[anchor=south west,inner sep=0] at (0.8,0.03) {\small{low}};
		\node[anchor=south west,inner sep=0] at (0.8,0.45) {\small{high}};
		\end{scope}
		\end{tikzpicture}
		}
	\hfill
		\subfigure[Temperature]{\includegraphics[trim=13cm 0mm 13.5cm 0.9cm,clip,width=0.218\textwidth]{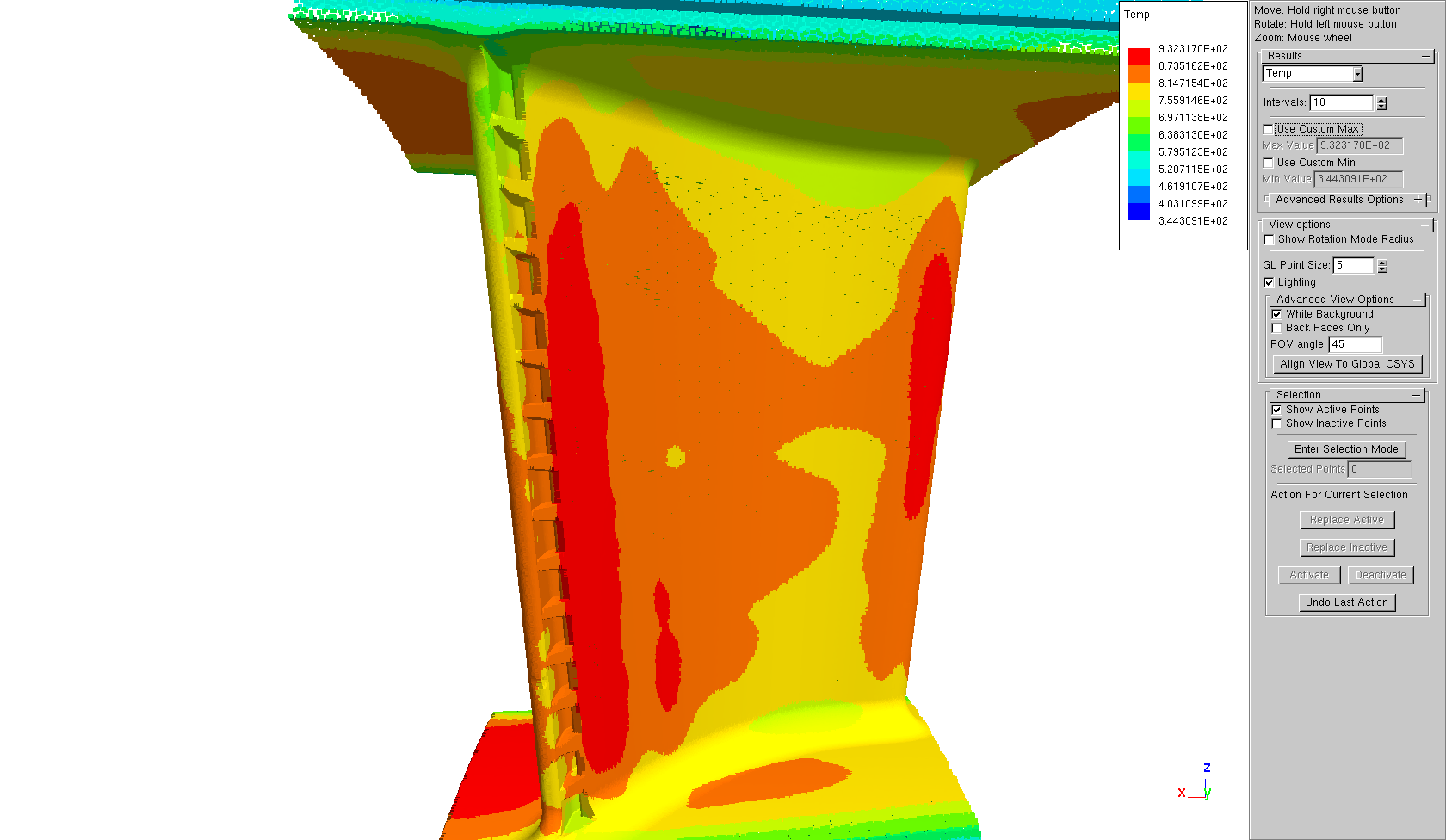}}
		\caption{NODAL FEA-VALUES OF TEMPERATURE AND VON MISES STRESS AT PRESSURE SIDE OF TURBINE VANE}
		\label{fig:TEMP_EQSTRESS}
		\vspace{-1\baselineskip}
\end{figure}

In this section the probabilistic model is applied to a turbine vane made of polycrystalline cast RENE 80. All probabilistic analyses of its LCF crack initiation life are based on FEA simulating a thermomechanical load as in the operating state. The model consists of tetrahedral elements and heat transfer analysis is performed assuming an undamaged system of thermal barrier coating (TBC) and bond coat (BC). Note that the present examination neglects the complex mechanical interaction between coating and substrate material. The structural analysis delivers the elastic strain tensor and temperature field at all nodes. From stress tensor data the elastic-plastic strain field $\epsilon(\mathbf{x})$ (Fig.~\ref{fig:TEMP_EQSTRESS} (a)) is calculated which is then further utilized for computing the deterministic life field $N_{i_\mathrm{det}}(\mathbf{x})$ according to Eqn.~\eqref{eq:locPM_locCMB} together with the temperature data.

Note that instead of nodal values all results of the FEA postprocessor, explicitly the surface integration in Eq.~\eqref{eq:locPM_ScaleSurfInt}, are evaluated and plotted at the coordinates of integration points of the finite elements\footnote{Quadratures of higher order than two are chosen to rule out numerical nonlinearities, compare \cite{ASME2013Paper}} at the domain surface $\partial\Omega$ by applying the corresponding shape functions for interpolation.

\subsection{Probabilistic Lifing Without Notch Support}\label{sec:TBV_noNSP}
\begin{figure}[htbp]
	\centering
	\begin{tikzpicture}
		\node[anchor=south west,inner sep=0] (haz_map) at (0,0) {\includegraphics[trim=0cm 0cm 0cm 2cm,clip,width=0.45\textwidth]{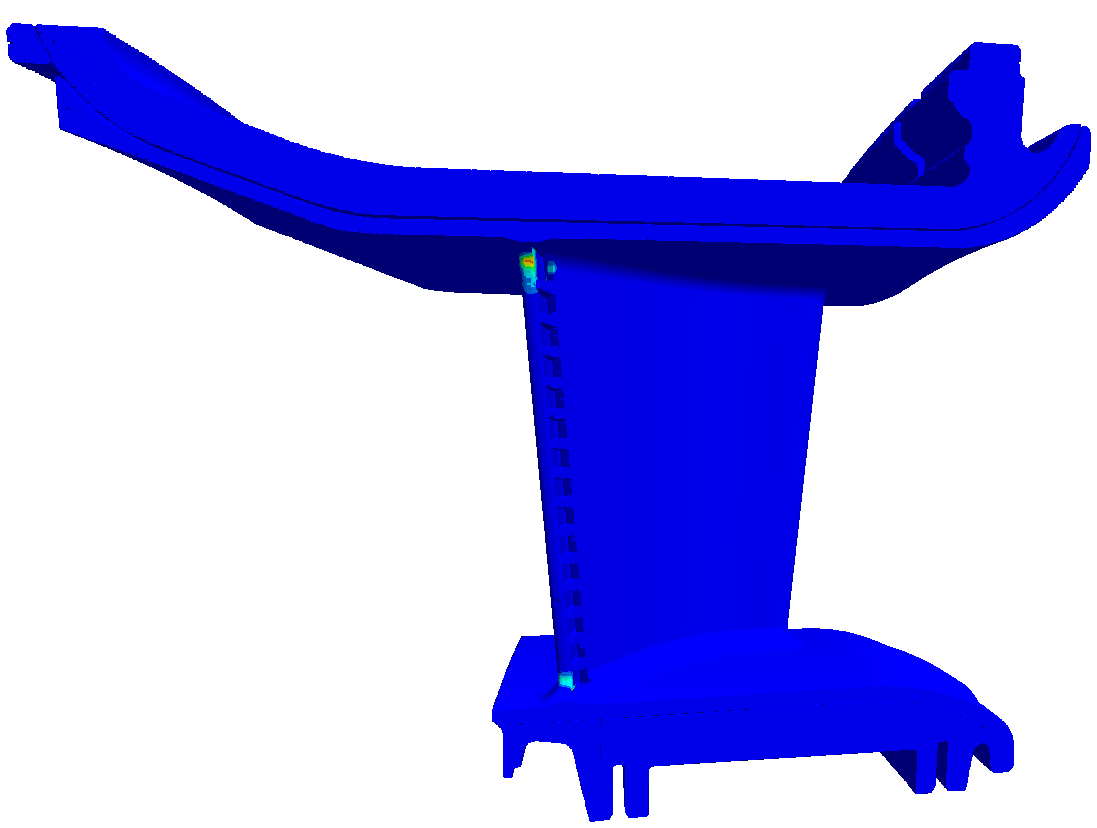}};
		\begin{scope}[x={(haz_map.south east)},y={(haz_map.north west)}]
		\draw[red,thick] (0.45,0.68) rectangle (0.55,0.78);
		\node[anchor=south west,inner sep=0] at (0,0.5) {\includegraphics[trim=25cm 20cm 19cm 7.5cm,clip,width=0.15\textwidth]{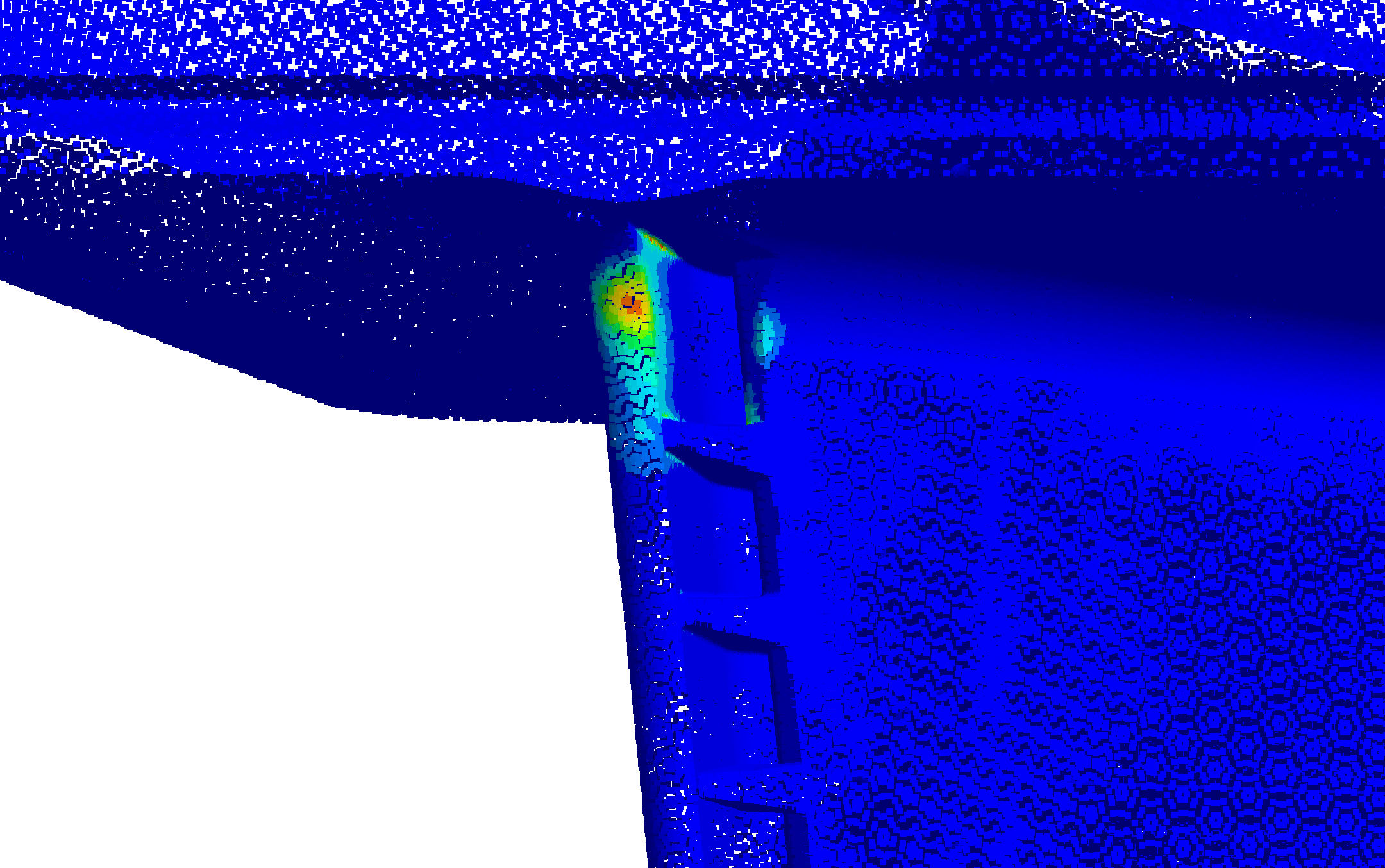}};
		\draw[red,thick] (0,0.5) rectangle (0.3334,0.8);
		\draw[red,thick] (0.45,0.78) -- (0.3334,0.8);
		\draw[red,thick] (0.45,0.68) -- (0.3334,0.5);
		\draw[red,thick] (0.48,0.15) rectangle (0.58,0.25);
		\node[anchor=south west,inner sep=0] at (0,0) {\includegraphics[trim=30cm 20cm 24cm 13.5cm,clip,width=0.15\textwidth]{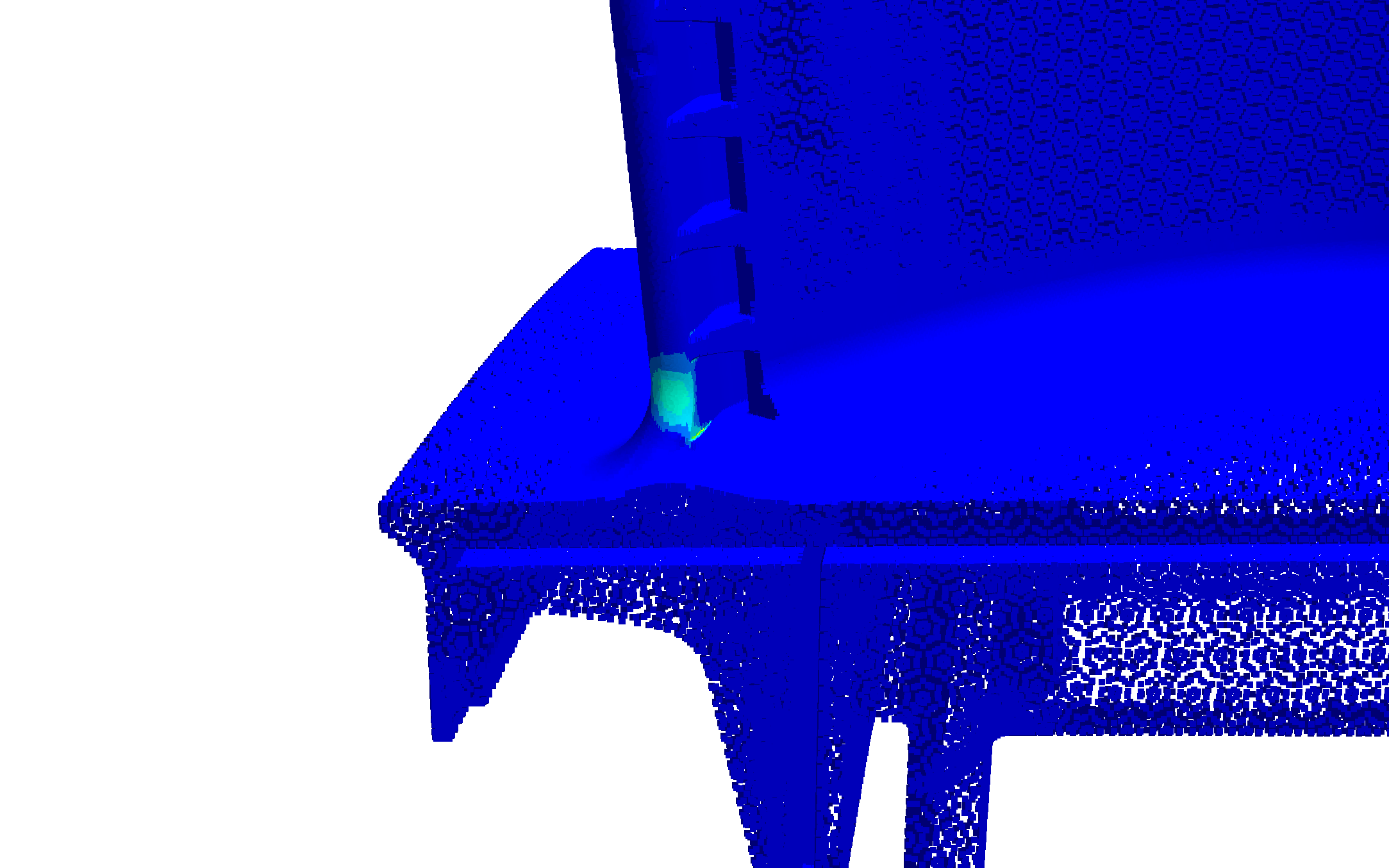}};
		\draw[red,thick] (0,0) rectangle (0.3334,0.3);
		\draw[red,thick] (0.48,0.25) -- (0.3334,0.3);
		\draw[red,thick] (0.48,0.15) -- (0.3334,0);
		\node[anchor=south west,inner sep=0] at (0.94,0.35) {\includegraphics[trim=77.5cm 35cm 12cm 2cm,clip,width=0.007\textwidth]{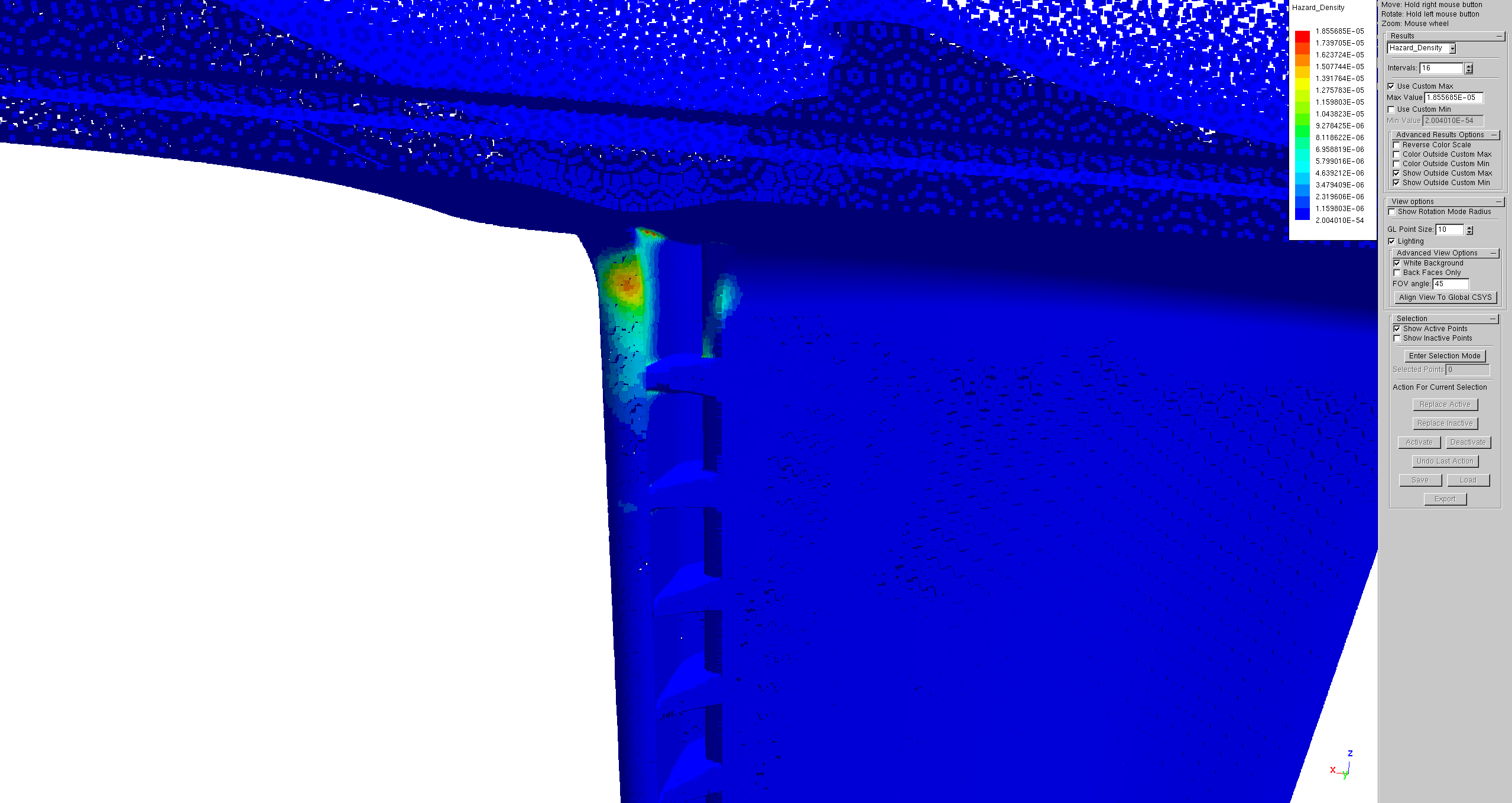}};
		\node[anchor=south west,inner sep=0] at (0.85,0.36) {\small{low}};
		\node[anchor=south west,inner sep=0] at (0.85,0.62) {\small{high}};
		\end{scope}
	\end{tikzpicture}
	\caption{HAZARD DENSITY PLOT OF VANE}
	\label{fig:HazDen_whole}
	\vspace{-1\baselineskip}
\end{figure}
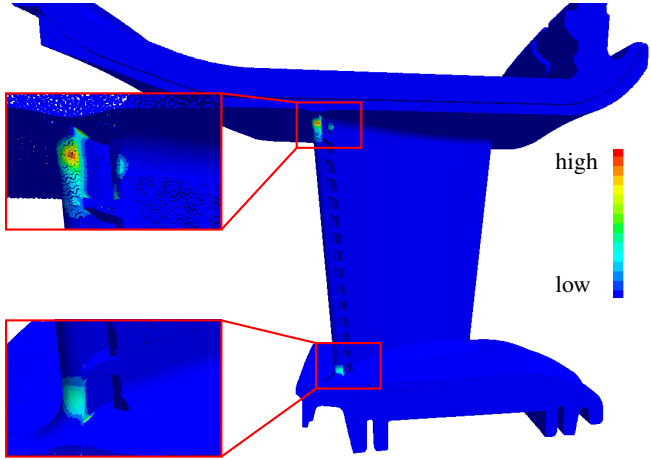

In this subsection the same crack initiation prediction model as presented in \cite{ASME2013Paper} is applied. Having calculated $N_{i_\mathrm{det}}(\mathbf{x})$, one can project the hazard density $(1/N_{i_\mathrm{det}}(\mathbf{x}))^m$ onto the vane as shown in Fig.~\ref{fig:HazDen_whole}.

It indicates areas of increased risk at the transitions from airfoil to inner shroud and outer shroud which occur only at the trailing edge. This strongly correlates with the locations of highest stress observed in Fig.~\ref{fig:TEMP_EQSTRESS} (a). However, since the material parameters are also temperature dependent, the vane temperature field (Fig.~\ref{fig:TEMP_EQSTRESS} (b)) influences the probability of crack initiation as well. 

From the hazard density and Weibull shape $m$ one can derive the probability distribution function for crack initiation events using Eqn.~\eqref{eq:DistrFunct} which is shown in Fig.~\ref{fig:cdf_both}
The ratio of probabilistic average life for LCF crack initiation and the deterministic life of a certain smooth specimen subjected to equal maximum strain $\epsilon_a$ is defined as the \textit{size effect factor}. This turbine vane examination results in a size effect factor of 3.25. Although the entire vane surface is much larger than the gauge section area of those LCF-specimens, the regions with critical hazard density at the vane are confined to very narrow spots. Thereby the overall hazard rate for the critical surface at the vane is smaller than the hazard rate for a standard LCF specimen and thus leads to higher probabilistic average life which corresponds to the statistical size effect.

\subsection{Probabilistic Lifing With Notch Support}\label{sec:TBV_NSP}
\begin{figure}[htbp]
	\centering
	\begin{tikzpicture}
		\node[anchor=south west,inner sep=0] (chi_map) at (0,0) {\includegraphics[trim=13cm 0mm 14cm 2cm,clip,width=0.4\textwidth]{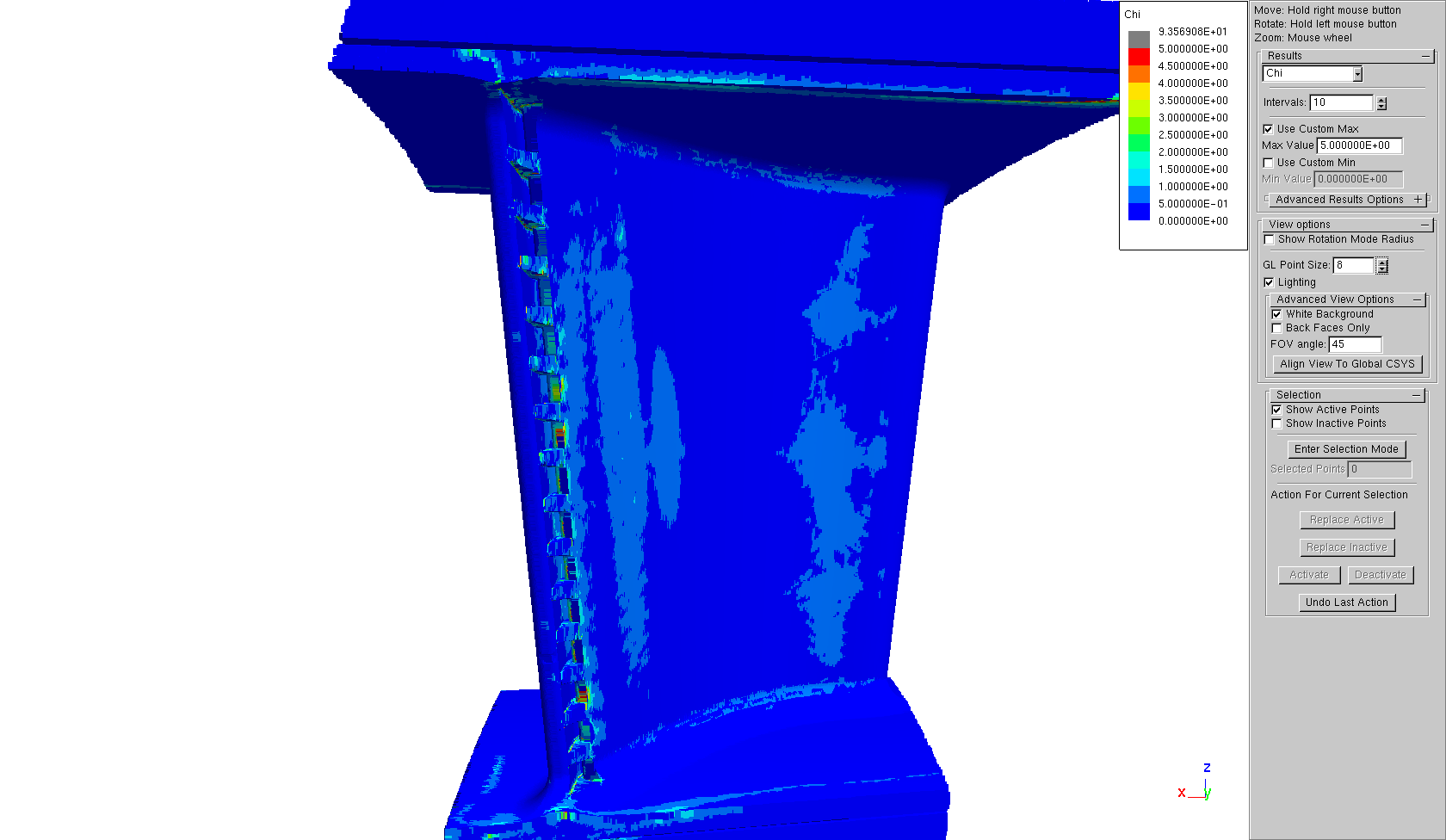}};
		\begin{scope}[x={(chi_map.south east)},y={(chi_map.north west)}]
		\node[anchor=south west,inner sep=0] at (0.95,0) {\includegraphics[trim=77.5cm 35.6cm 12cm 1.9cm,clip,width=0.01\textwidth]{Temp_Zoom_PS10.png}};
		\node[anchor=south west,inner sep=0] at (0.85,0.02) {\small{low}};
		\node[anchor=south west,inner sep=0] at (0.85,0.38) {\small{high}};
		\end{scope}
		\end{tikzpicture}
	\caption{$\chi$ VALUES AT THE PRESSURE SIDE OF AIRFOIL, LINERARILY SCALED}
	\label{fig:Chi_plot}
	\vspace{-1\baselineskip}
	\end{figure}

In order to consider the notch support effect for the probabilistic crack initiation life, the lifing algorithm calculates the analytical derivatives of the FEA shape functions to obtain the von Mises stress gradient at the domain surface as in Eqn.~\eqref{eq:NSE_Chicalc}. The resulting $\chi$ field at the pressure side of the airfoil is shown in Fig.~\ref{fig:Chi_plot}.

High values of the related stress gradient occur at sharp shape transitions of the geometry, for example at the crosspieces of the cooling channel outlets and the fillet radii. Negative $\chi$-values (strain decreasing towards surface) are set to zero.
	
\begin{figure}[htbp]
	\centering
	\includegraphics[trim=0cm 0cm 0cm 0cm,clip,width=0.25\textwidth]{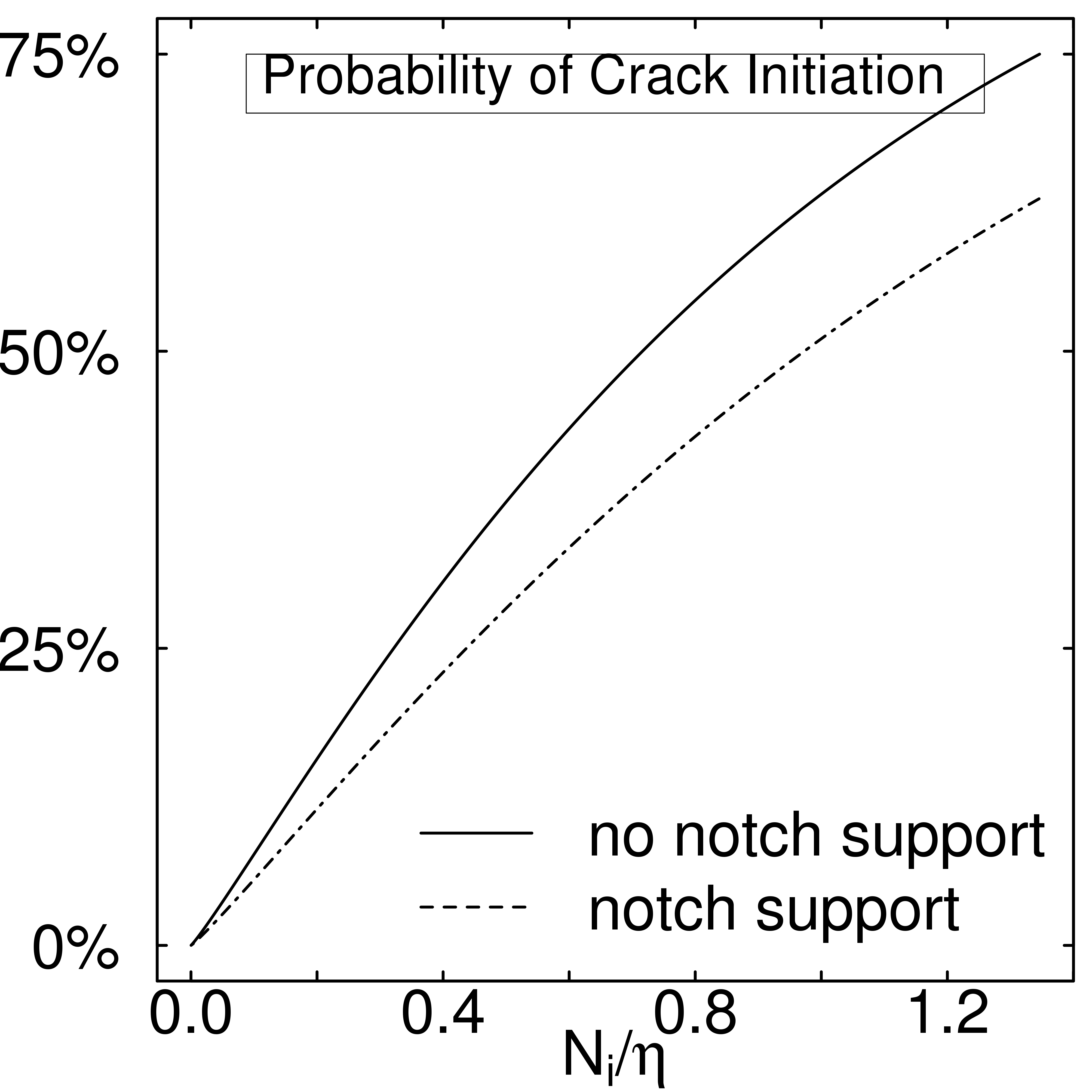}
	\caption{PROBABILITIES OF CRACK INITIATION AT VANE NEGLECTING AND CONSIDERING NOTCH SUPPORT}
	\label{fig:cdf_both}
\end{figure}

If the notch support effect is considered in the calculation, the scale value $\eta$ changes according to Eqn.~\eqref{eq:locPM_ScaleSurfInt} and therefore shifts the probability distribution for crack initiation. The difference in the distribution functions for crack initiation for the whole vane, computed considering and neglecting notch support, can be seen in Fig.~\ref{fig:cdf_both}. When enabling the notch support in the local probabilistic model by using Eqn.~\eqref{eq:locPM_locCMB_andNS} instead of Eqn.~\eqref{eq:locPM_locCMB}, a larger Weibull scale parameter for the distribution is received. This reduces the slope of the distribution function and results in \unit{36}{\%} higher life.

\begin{figure}[htbp]
	\subfigure[no notch support]{
	\begin{tikzpicture}
		\node[anchor=south west,inner sep=0] (Notch) at (0,0) {\includegraphics[trim=20cm 8cm 25cm 9cm,clip,width=0.23\textwidth]{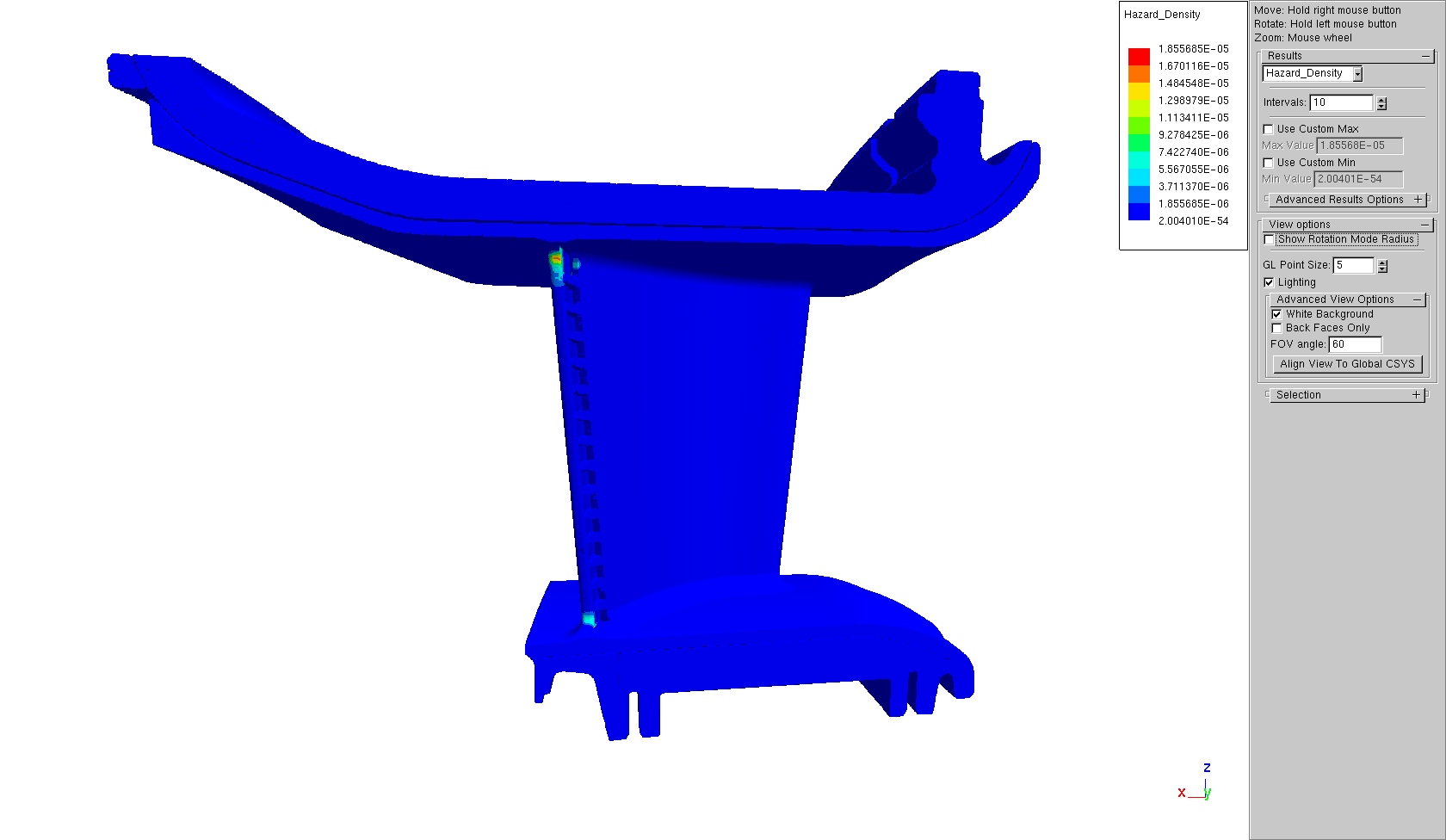}};
	\begin{scope}[x={(Notch.south east)},y={(Notch.north west)}]
	\node[anchor=south west,inner sep=0] at (0.15,0.2) {\includegraphics[trim=77.5cm 35.2cm 12cm 2cm,clip,width=0.007\textwidth]{Basic_Hazard_Zoom_PS10.png}};
		\node[anchor=south west,inner sep=0] at (0,0.21) {\small{low}};
		\node[anchor=south west,inner sep=0] at (0,0.5) {\small{high}};
		\end{scope}
	\end{tikzpicture}
	}
	\subfigure[notch support]{\includegraphics[trim=20cm 8cm 25cm 9cm,clip,width=0.23\textwidth]{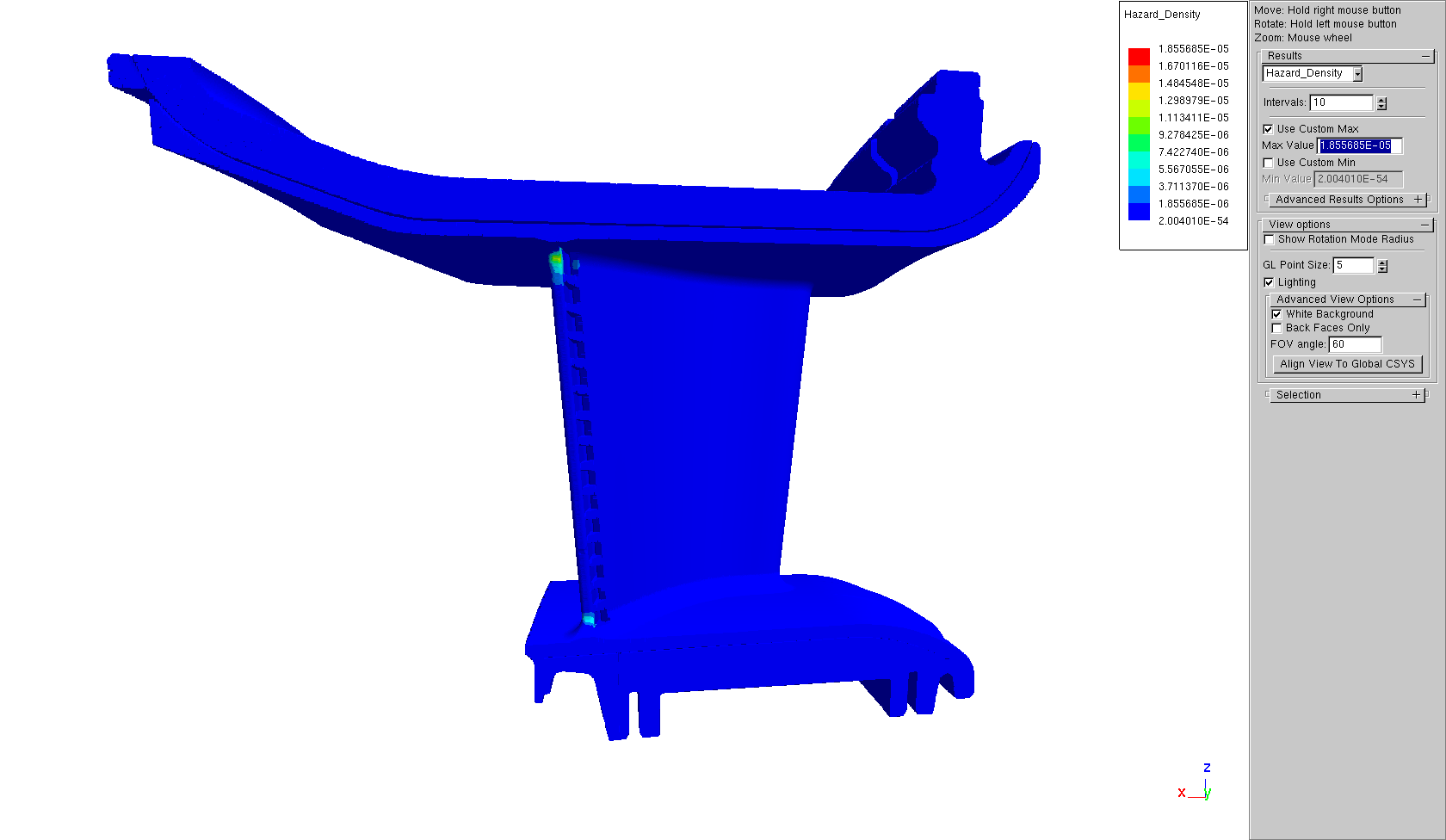}}
	\caption{HAZARD DENSITIES AT AIRFOIL FROM PREDICTIONS NEGLECTING AND CONSIDERING NOTCH SUPPORT}
	\label{fig:HazDen_both}
\end{figure}

\begin{figure}[htbp]
\subfigure[Strain]{
	\begin{tikzpicture}
		\node[anchor=south west,inner sep=0] (strainmap) at (0,0) {\includegraphics[trim=20cm 20cm 45cm 5cm,clip,width=0.23\textwidth]{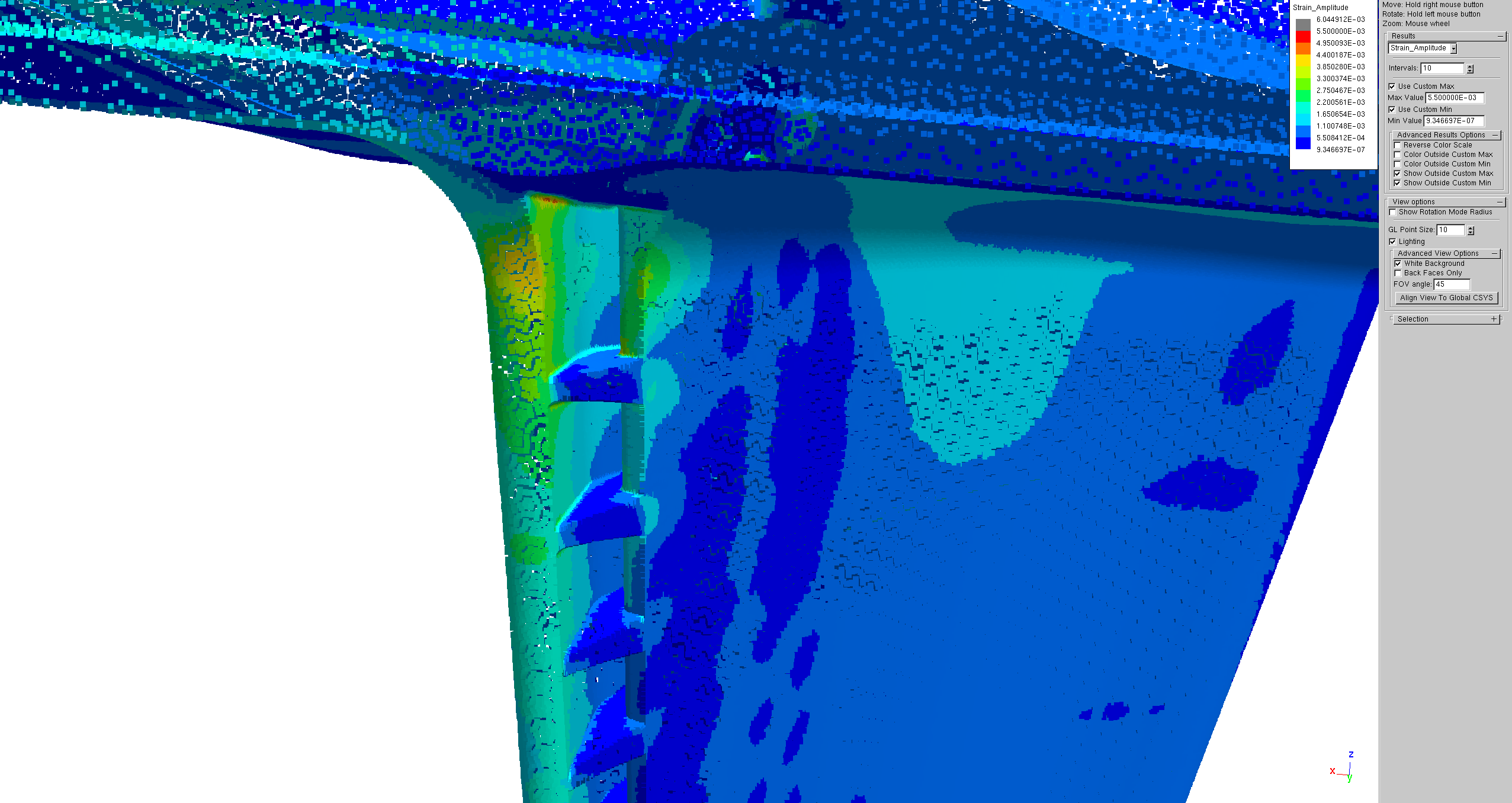}};
		\begin{scope}[x={(strainmap.south east)},y={(strainmap.north west)}]
		\draw[red,thick,rounded corners] (0.36,0.39) rectangle (0.52,0.61);
		\node[red] at (0.63,0.55) {spot 1};
		\draw[red,thick,rounded corners] (0.45,0.66) rectangle (0.58,0.74);
		\node[red] at (0.7,0.7) {spot 2};
		\node[anchor=south west,inner sep=0] at (0,0) {\includegraphics[trim=77.5cm 35cm 12cm 1cm,clip,width=0.007\textwidth]{Temp_Zoom_PS10.png}};
		\node[anchor=south west,inner sep=0] at (0.05,0.05) {\small{low}};
		\node[anchor=south west,inner sep=0] at (0.05,0.52) {\small{high}};
		\end{scope}
	\end{tikzpicture}
	}
\subfigure[Temperature]{\includegraphics[trim=20cm 20cm 45cm 5cm,clip,width=0.23\textwidth]{Temp_Zoom_PS10.png}}
\caption{ELASTIC STRAIN AND TEMPERATURE FIELD IN TRANSITION FROM TRAILING EDGE OF AIRFOIL TO OUTER SHROUD}
\label{fig:Temp-Strain}
\vspace{-1\baselineskip}
\end{figure}

Comparing the hazard densities in Fig.~\ref{fig:HazDen_both}, one can notice the similar shape of the risk patches but at the same time decreased values in the results of the notch support examination. Closer examination of the upper section of the airfoil's trailing edge is given in Fig.~\ref{fig:Temp-Strain}. 

\begin{figure}[htbp]
	\centering
	\subfigure[von Mises Stress]{
		\begin{tikzpicture}
		\node[anchor=south west,inner sep=0] (strainmap) at (0,0) {\includegraphics[trim=26cm 25cm 45cm 8cm,clip,width=0.23\textwidth]{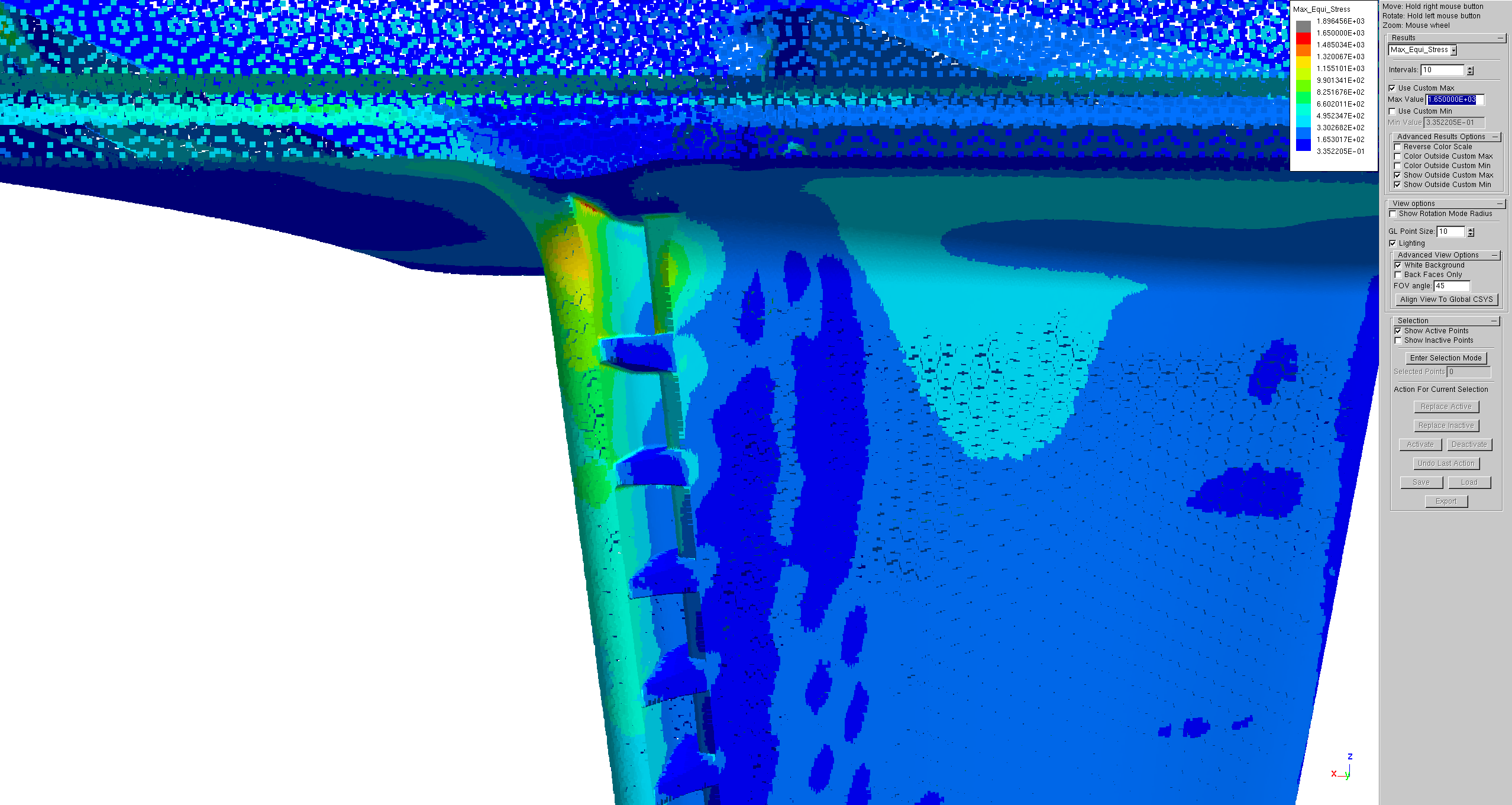}};
		\begin{scope}[x={(strainmap.south east)},y={(strainmap.north west)}]
		\draw[red,thick,rounded corners,rotate around={8:(0.475,0.45)}] (0.35,0.27) rectangle (0.5,0.63);
		\node[red] at (0.62,0.55) {spot 1};
		\draw[red,thick,rounded corners,rotate around={-10:(0.475,0.7)}] (0.4,0.65) rectangle (0.55,0.75);
		\node[red] at (0.7,0.7) {spot 2};
		\end{scope}
	\end{tikzpicture}
	}
	\subfigure[$\chi$]{\includegraphics[trim=29.5cm 29cm 46.5cm 8cm,clip,width=0.23\textwidth]{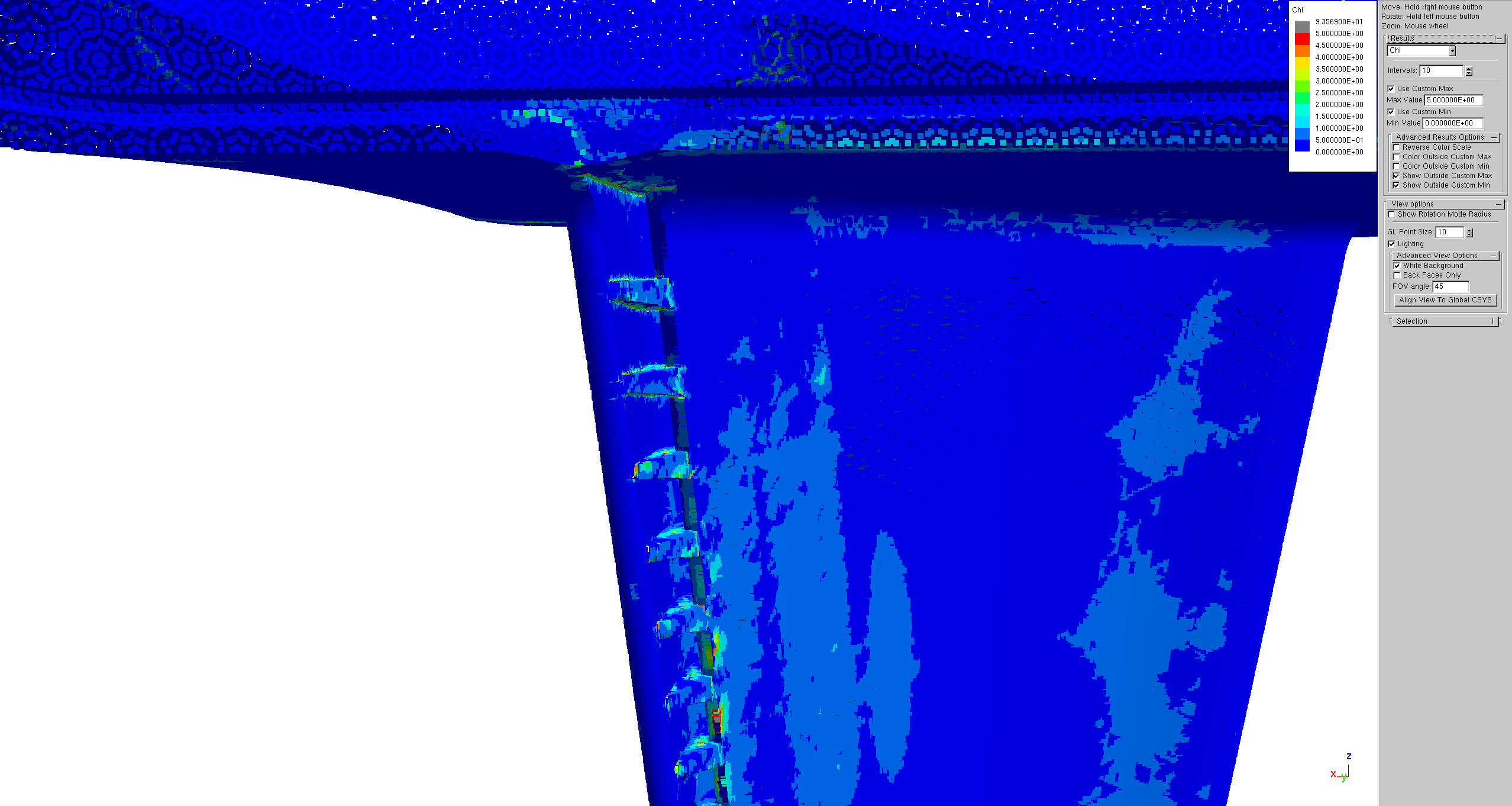}}
	\subfigure[$\chi_T$]{
		\begin{tikzpicture}
		\node[anchor=south west,inner sep=0] (chiT) at (0,0) {\includegraphics[trim=26cm 29cm 46cm 5cm,clip,width=0.23\textwidth]{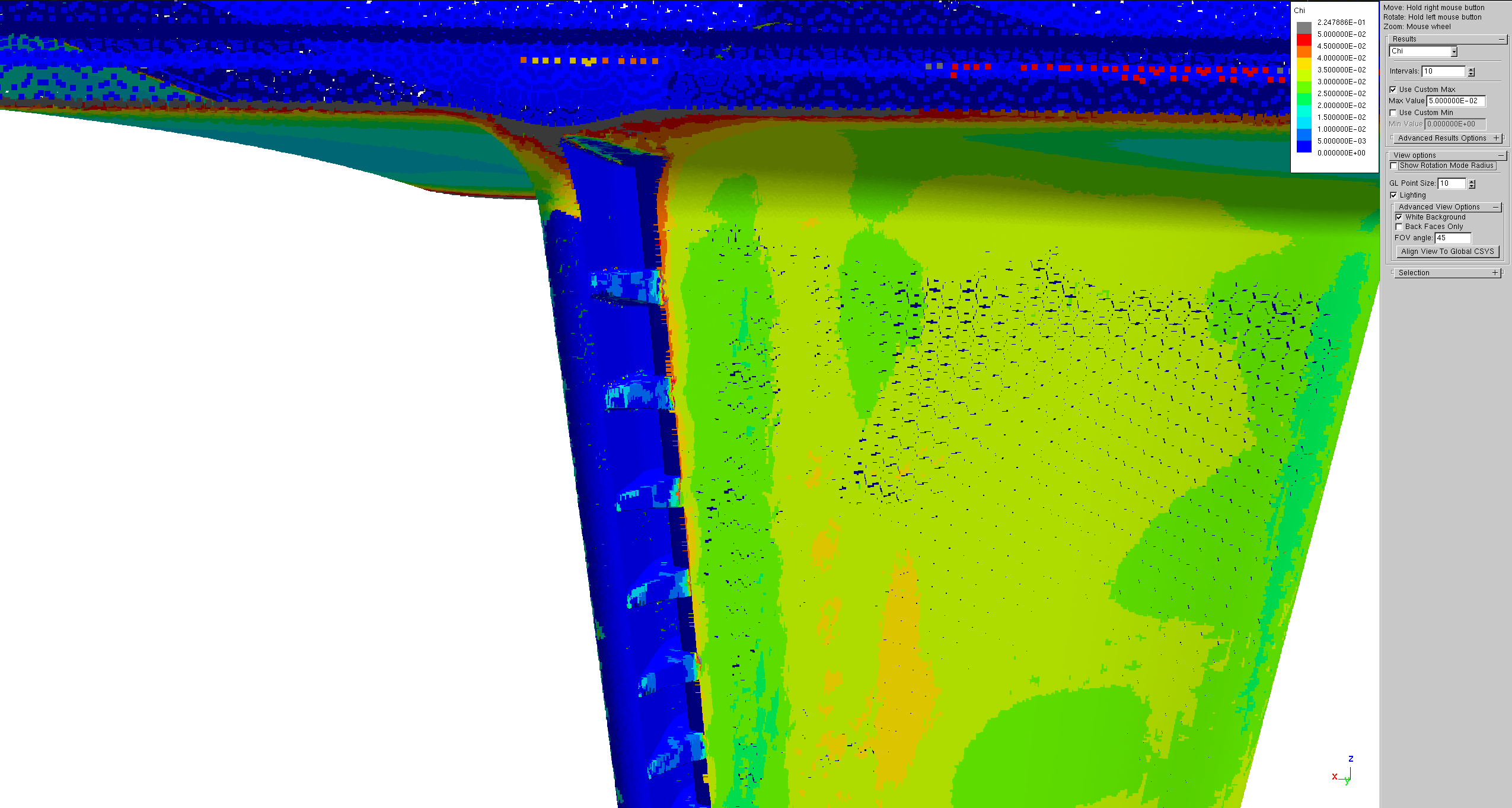}};
		\begin{scope}[x={(chiT.south east)},y={(chiT.north west)}]
		\node[anchor=south west,inner sep=0] at (0,0) {\includegraphics[trim=77.5cm 35cm 12cm 1.2cm,clip,width=0.007\textwidth]{Temp_Zoom_PS10.png}};
		\node[anchor=south west,inner sep=0] at (0.05,0.05) {\small{low}};
		\node[anchor=south west,inner sep=0] at (0.05,0.6) {\small{high}};
		\end{scope}
		\end{tikzpicture}
	}
	\caption{VON MISES STRESS, $\chi$ AND $\chi_T$ FIELD IN TRANSITION FROM TRAILING EDGE OF AIRFOIL TO OUTER SHROUD}
	\label{fig:Stress-Chi}
	\vspace{-1\baselineskip}
\end{figure}

Fig.~\ref{fig:Temp-Strain} (a) shows two spots of distinctively visible strain concentrations originating mostly from inhomogeneous thermal expansion. This is illustrated in Fig.~\ref{fig:Temp-Strain} (b) and Fig.~\ref{fig:Stress-Chi} (c) which show the normalized temperature gradient $\chi_T$. The local strains, seen in Fig.~\ref{fig:Temp-Strain} (b), cause the stresses in the respective locations, as shown in Fig.~\ref{fig:Stress-Chi} (a). 

The related stress gradient $\chi$ is mapped onto the geometry of the examined section in Fig.~\ref{fig:Stress-Chi} (b). However the distinctive larger spot of high stress at the trailing edge near the outer shroud (spot 1) is not recognizable in the $\chi$-field. Since the shape in that area is relatively smooth compared to the spot in the edge of the top cooling channel, the local $\chi$-values are not significantly elevated compared to the surrounding material. In contrast, higher $\chi$-values occur in spot 2 since it features a sharp geometry transition causing a very inhomogeneous stress distribution. The notch support effect should therefore have a higher impact on the probabilistic crack initiation life in spot 2 compared to spot 1.

\begin{figure}[htbp]
	\centering
	\vspace{-0.5\baselineskip}
	\subfigure[no notch support]{
		\begin{tikzpicture}
		\node[anchor=south west,inner sep=0] (haz_map1) at (0,0) {\includegraphics[trim=8cm 18cm 38cm 4cm,clip,width=0.22\textwidth]{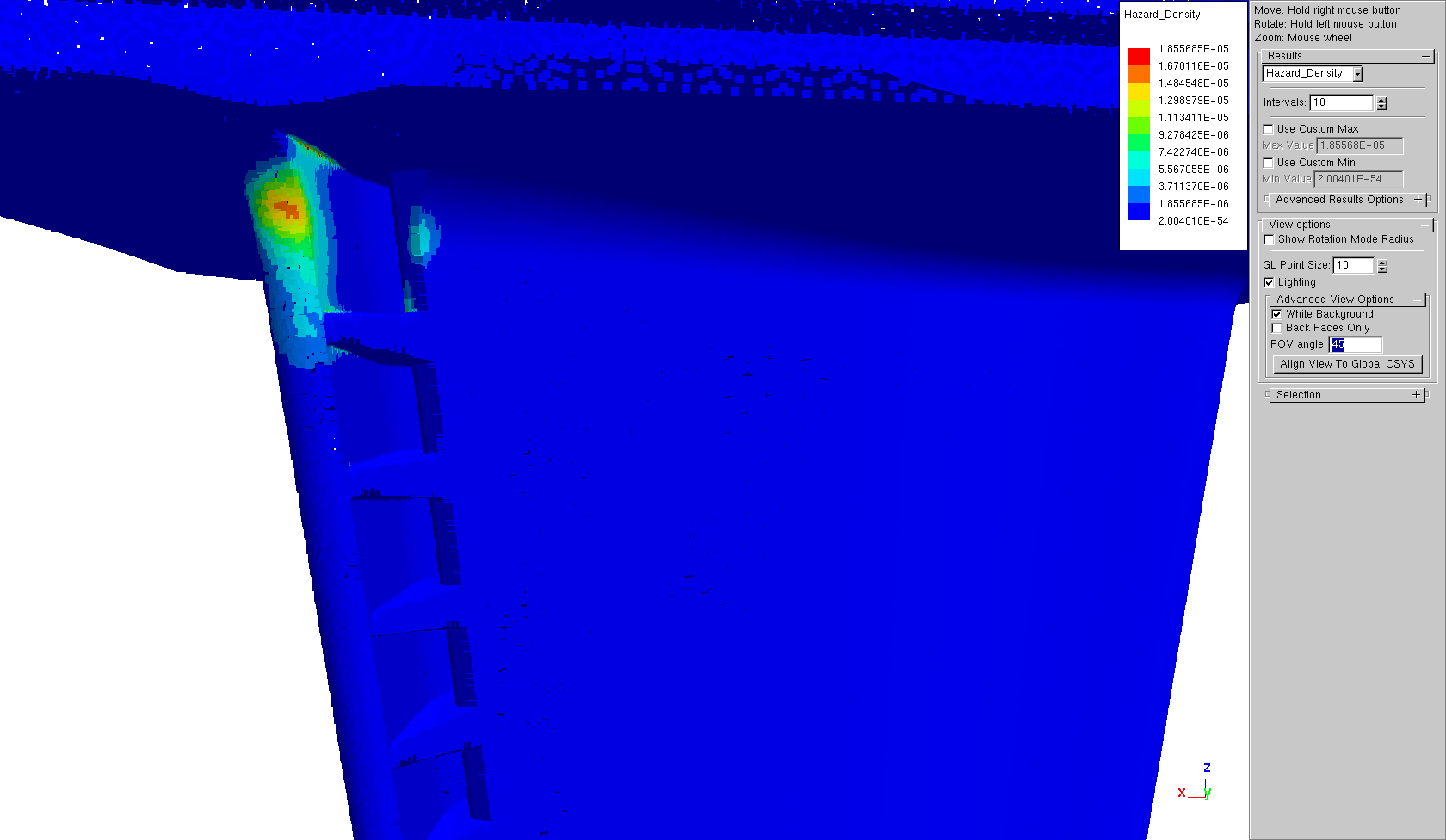}};
		\begin{scope}[x={(haz_map1.south east)},y={(haz_map1.north west)}]
		\draw[red,thick,rounded corners,rotate around={10:(0.29,0.56)}] (0.18,0.44) rectangle (0.4,0.78);
		\node[red] at (0.13,0.8) {spot 1};
		\draw[red,thick,rounded corners,rotate around={-30:(0.365,0.825)}] (0.25,0.79) rectangle (0.48,0.86);
		\node[red] at (0.45,0.92) {spot 2};
		\node[anchor=south west,inner sep=0] at (0.54,0.3) {\includegraphics[trim=13cm 4cm 11cm 10cm,clip,width=0.1\textwidth]{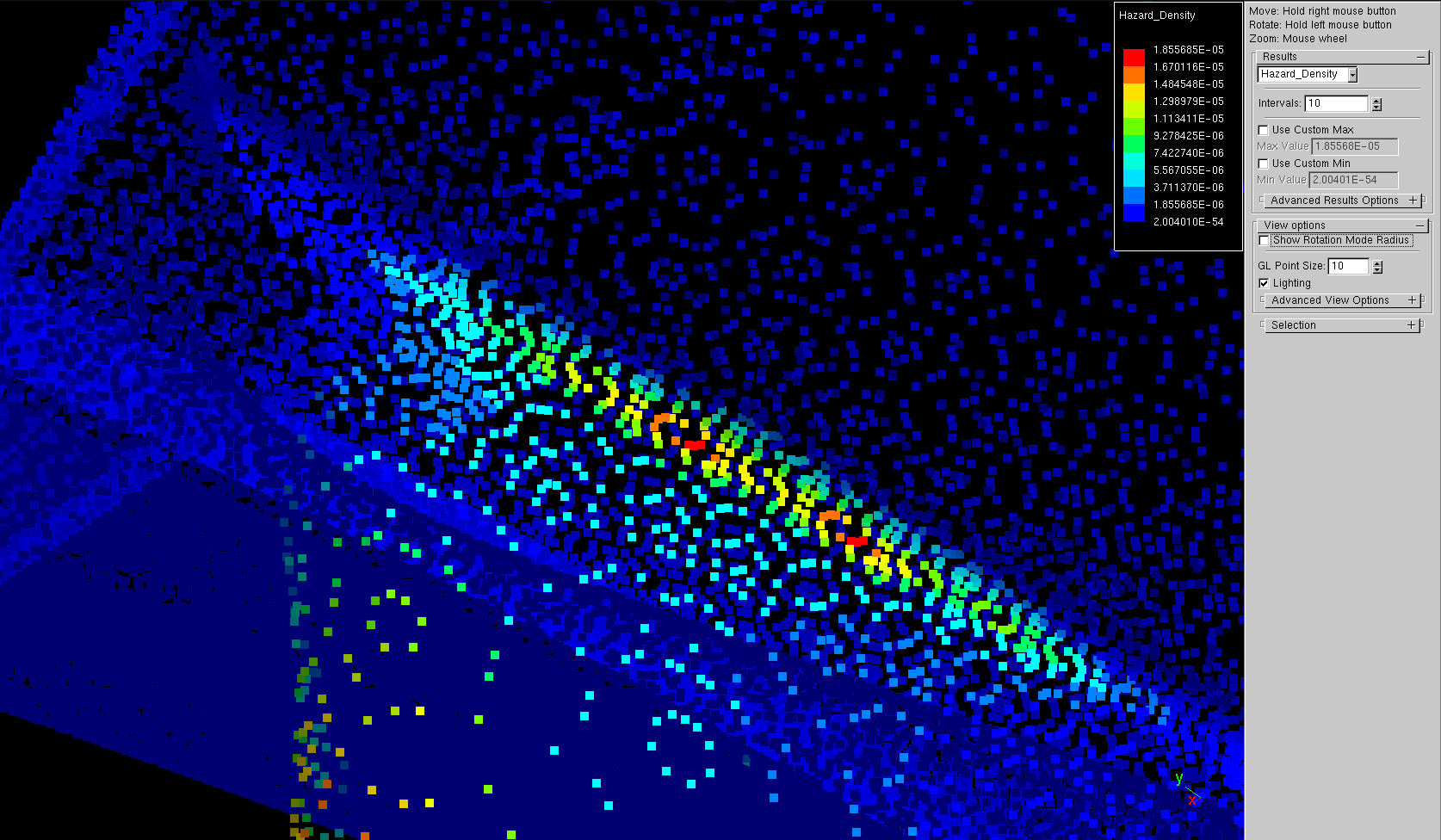}};
		\draw[red,thick] (0.54,0.3) rectangle (1,0.59);
		\draw[red,thick] (0.465,0.79) -- (1,0.59);
		\draw[red,thick] (0.45,0.75) -- (0.54,0.3);
		\end{scope}
	\end{tikzpicture}
	}
	\subfigure[notch support]{
	\begin{tikzpicture}
		\node[anchor=south west,inner sep=0] (haz_map2) at (0,0) {\includegraphics[trim=8.2cm 15.3cm 34.8cm 3.8cm,clip,width=0.22\textwidth]{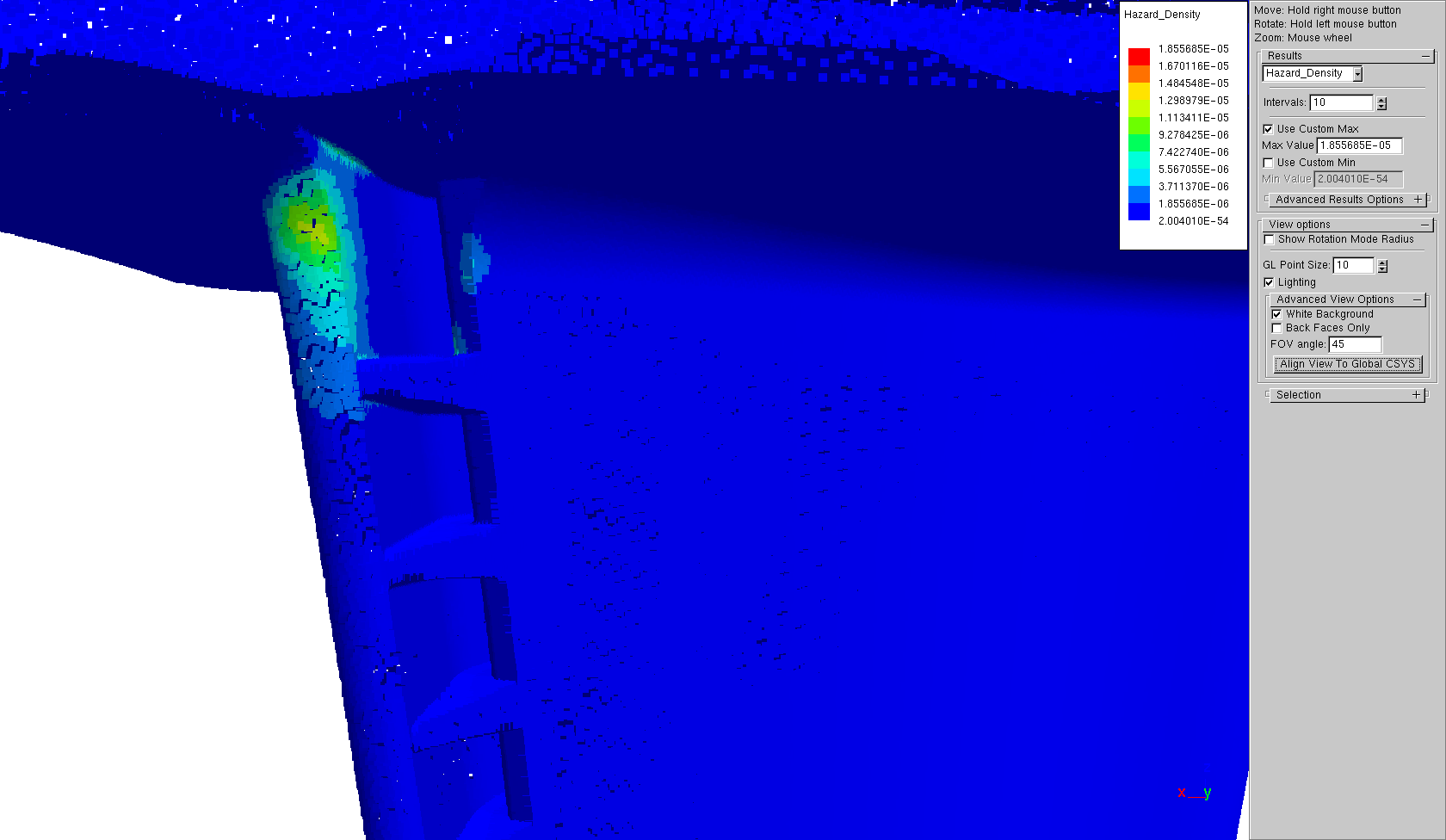}};
		\begin{scope}[x={(haz_map2.south east)},y={(haz_map2.north west)}]
		\draw[red,thick,rounded corners,rotate around={10:(0.29,0.56)}] (0.18,0.44) rectangle (0.4,0.78);
		\node[red] at (0.13,0.8) {spot 1};
		\draw[red,thick,rounded corners,rotate around={-30:(0.365,0.825)}] (0.25,0.79) rectangle (0.48,0.86);
		\node[red] at (0.45,0.92) {spot 2};
		\node[anchor=south west,inner sep=0] at (0.54,0.3) {\includegraphics[trim=13cm 9cm 20cm 10cm,clip,width=0.1\textwidth]{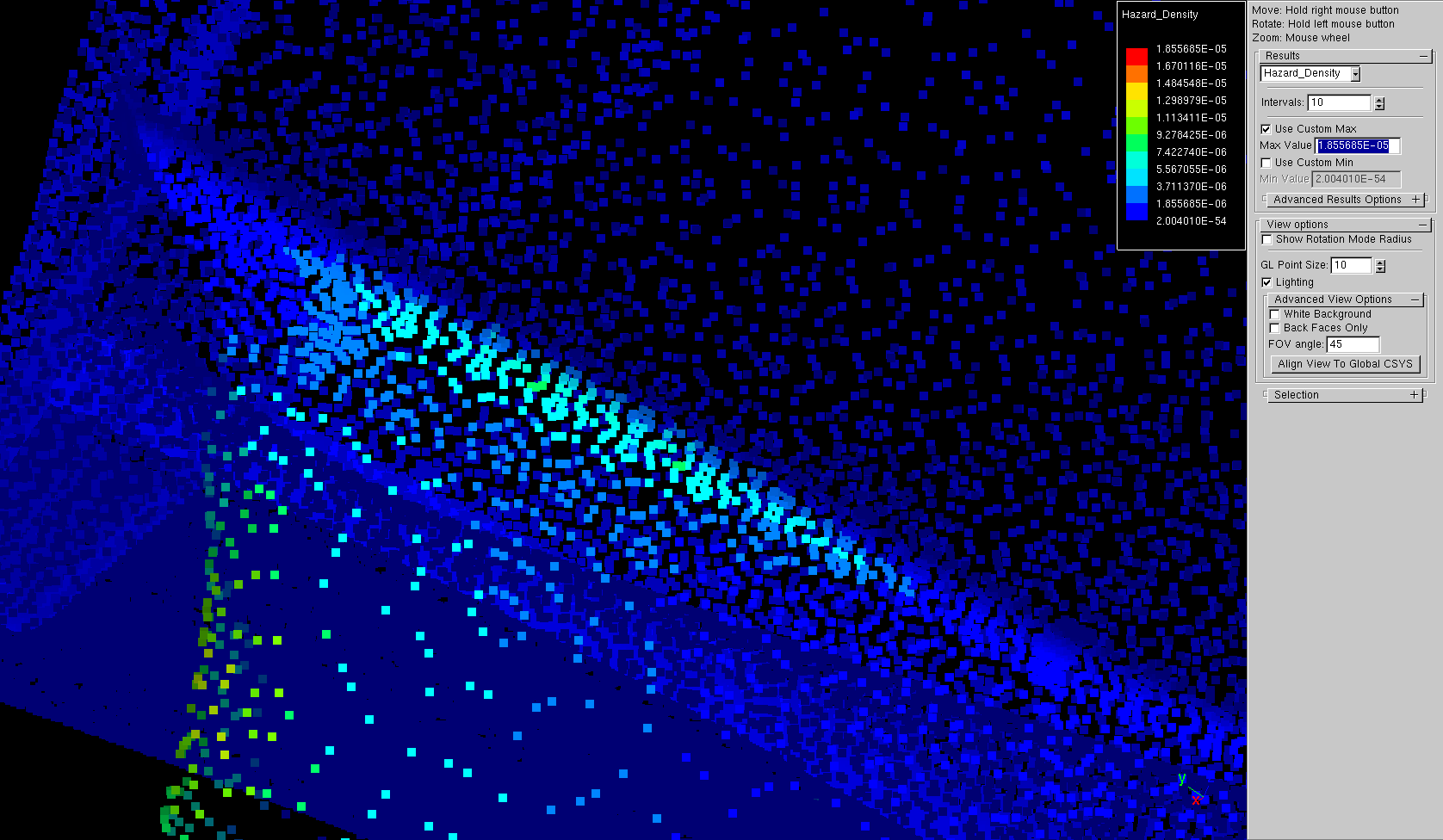}};
		\draw[red,thick] (0.54,0.3) rectangle (1,0.59);
		\draw[red,thick] (0.465,0.79) -- (1,0.59);
		\draw[red,thick] (0.45,0.75) -- (0.54,0.3);
		\node[anchor=south west,inner sep=0] at (0,0) {\includegraphics[trim=77.5cm 35cm 12.1cm 2cm,clip,width=0.008\textwidth,height=0.085\textheight]{Basic_Hazard_Zoom_PS10.png}};
		\node[anchor=south west,inner sep=0] at (0.05,0) {\small{low}};
		\node[anchor=south west,inner sep=0] at (0.05,0.38) {\small{high}};
		\end{scope}
	\end{tikzpicture}	
	}
	\caption{HAZARD DENSITY COMPARISON IN TRANSITION FROM TRAILING EDGE OF AIRFOIL TO OUTER SHROUD}
	\label{fig:Zoom_HazardDens}
	\vspace{-0.5\baselineskip}
\end{figure}

\begin{figure}[htbp]
	\centering
		\subfigure[Spot 1]{\includegraphics[trim=0cm 0cm 0cm 0cm,clip,width=0.23\textwidth]{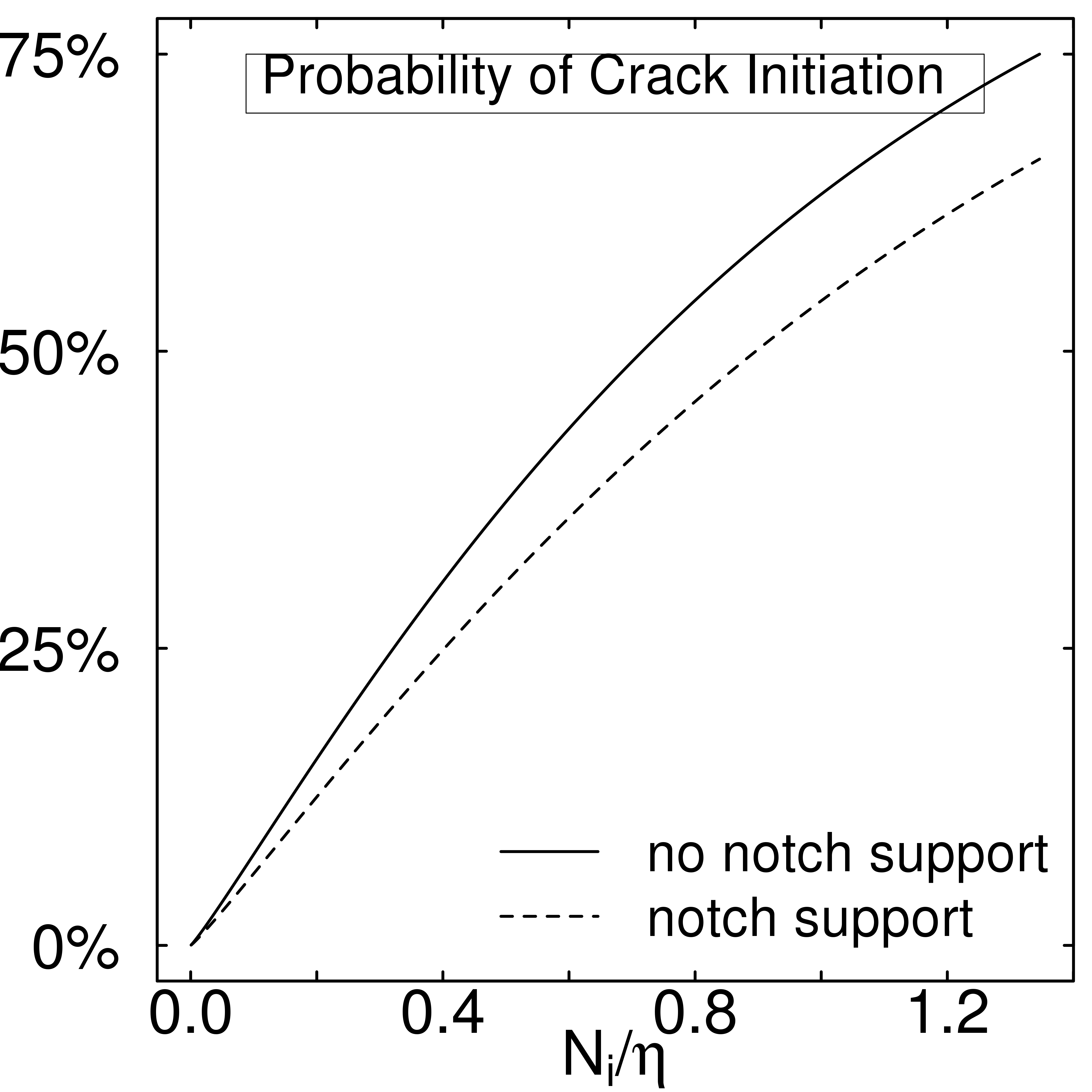}}
		\subfigure[Spot 2]{\includegraphics[trim=0cm 0cm 0cm 0cm,clip,width=0.23\textwidth]{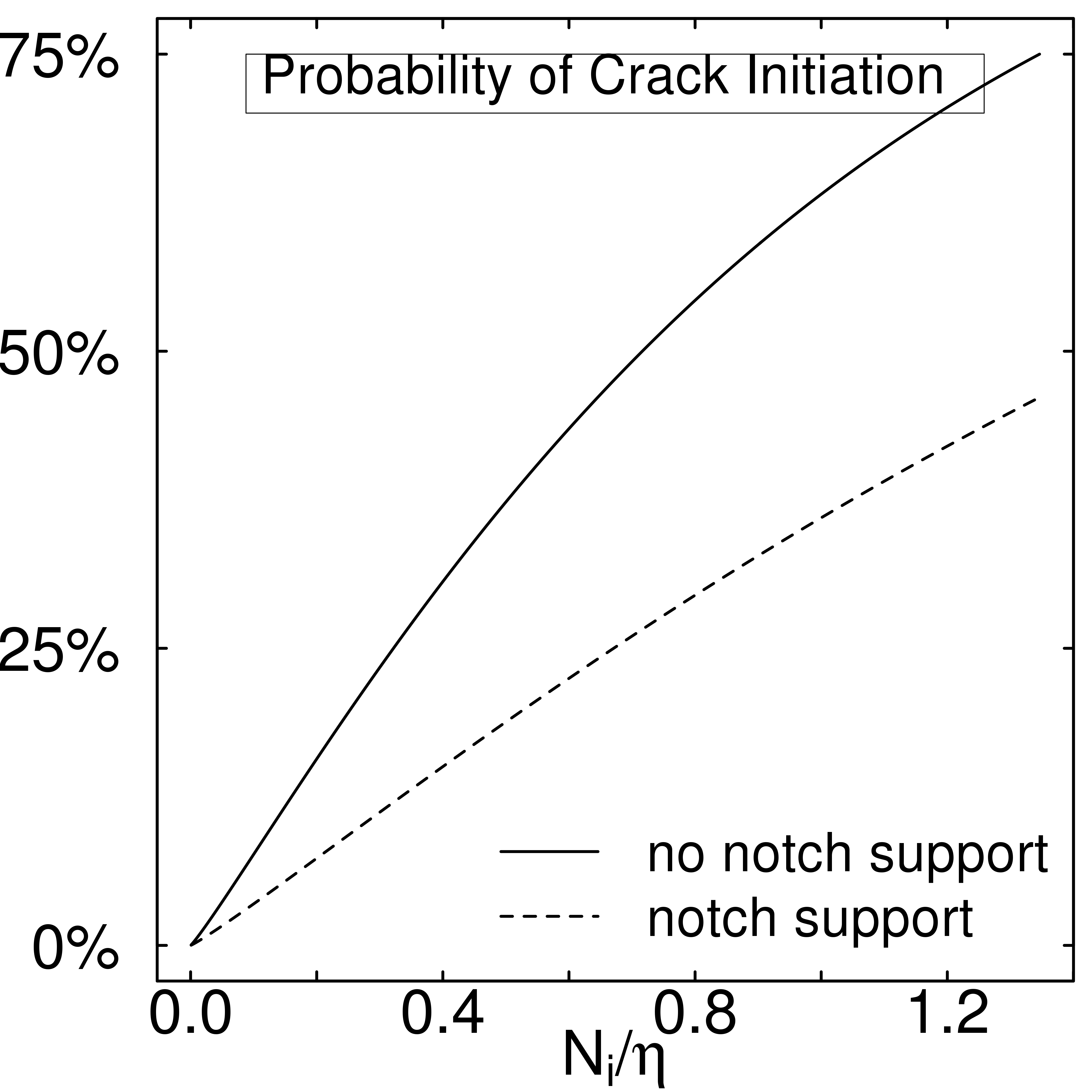}}
	\caption{COMPARISON OF RISK DECREASE IN CRITICAL SPOTS DUE TO NOTCH SUPPORT}
	\label{fig:PoCI_spot1u2}
	\vspace{-0.5\baselineskip}
\end{figure}

Indeed, the hazard densities in spot 2 are lower than in spot 1 of Fig.~\ref{fig:Zoom_HazardDens} (b). This is opposite to the situation in Fig.~\ref{fig:Zoom_HazardDens} (a) where higher hazard densities than in spot 1 are observed in spot 2. Computing the probabilities of crack initiation only from selected integration points in the respective spots confirms this observation. Fig.~\ref{fig:PoCI_spot1u2} (b) shows a significantly larger decrease in risk due to notch support for spot 2 compared to spot 1 in Fig.~\ref{fig:PoCI_spot1u2} (a).

The size effect factor in this crack initiation prediction, which is approximately 4.43, is now a \textit{combined size effect factor} because the notch support effect is incorporated in the hazard density and Weibull scale computation. Solving Eqn.~\eqref{eq:locPM_locCMB_andNS} for $N_{i_\mathrm{det}}(\mathbf{x})$ to use in the surface integration in Eq.~\eqref{eq:etaInt_HazDen} leads to higher values, i.e. less hazard density and higher probabilistic average life, compared to using Eqn.~\eqref{eq:locPM_locCMB}. The algorithm does not only consider stresses at the surface, but also the stress gradient which links to the stress field in the volume below the surface. The combined size effect is therefore increasing.

%% file: Conclusion.tex
\section*{DISCUSSION AND CONCLUSION}
In Section \ref{sec:LCF} the local probabilistic approach to LCF crack initiation prediction that was previously presented in \cite{Schmitz_Seibel} and \cite{ASME2013Paper} is reviewed. As specified in Subsection~\ref{sec:LCF_locPM}, the assumption of locally confined, independent LCF crack initiation events at engineering parts of polycrystalline metal allows a local hazard density approach. Using a Weibull distribution for the number of load cycles until crack initiation, the approach leads to a surface integral over the hazard density for the scale parameter which pays regard to the statistical size effect. Subsection~\ref{sec:LCF_NSE} shows how the notch support effect is incorporated in this model in a way that it combines with the statistical size effect. Calibration and validation of the presented notch support approach is exemplary outlined in Subsection~\ref{subsec:Model_Validation}. Both versions, the old and the extended, are applied for crack initiation life prediction of a gas turbine vane in Section~\ref{sec:LifePredictVane}. Areas of high hazard density on the vane are confined to small regions at the transitions of the trailing edge to inner and outer shroud and to edges in the cooling channel outlets. 

As expected, the model taking notch support into account predicts higher probabilistic life for the vane. By evaluating the critical spot at the trailing edge and the spot in the cooling channel edge separately, it is shown that significantly more probabilistic life is predicted for the second spot. That correlates with the philosophy of the implemented notch support model which states that larger stress gradients in a body lead to higher crack initiation life compared to a model that only considers homogeneous stress fields. Also, as shown in Fig.~\ref{fig:Fit_PredCurve}, the model is in good  agreement with experimental evidence.

However, the true LCF life of turbine components in operation is influenced by more factors such as TBC thickness and spallation, BC-substrate interaction behavior, thermomechanical fatigue (TMF) effects, grain size distribution, creep, manufacturing tolerances, variations in the operating conditions and the uncertainties of the parameter estimate $\mathbf{\hat{\theta}}$. The model applied for the work presented here is not specifically considering those. A probabilistic framework that combines the TBC/BC system life with the structural base material life and extends the present LCF-based model to TMF and even connects to fracture mechanics \cite{amann2016method} yet poses considerable future tasks for probabilistic gas turbine life prediction.